# A REGRESSION-BASED MONTE CARLO METHOD TO SOLVE BACKWARD STOCHASTIC DIFFERENTIAL EQUATIONS[1]

By Emmanuel Gobet, Jean-Philippe Lemor and Xavier Warin

*Centre de Mathématiques Appliquées, Electricité de France and Électricité de France*

We are concerned with the numerical resolution of backward stochastic differential equations. We propose a new numerical scheme based on iterative regressions on function bases, which coefficients are evaluated using Monte Carlo simulations. A full convergence analysis is derived. Numerical experiments about finance are included, in particular, concerning option pricing with differential interest rates.

**1. Introduction.** In this paper we are interested in numerically approximating the solution of a decoupled forward–backward stochastic differential equation (FBSDE)

$$S_t = S_0 + \int_0^t b(s, S_s)\,ds + \int_0^t \sigma(s, S_s)\,dW_s, \tag{1}$$

$$Y_t = \Phi(\mathbf{S}) + \int_t^T f(s, S_s, Y_s, Z_s)\,ds - \int_t^T Z_s\,dW_s. \tag{2}$$

In this representation, $\mathbf{S} = (S_t : 0 \leq t \leq T)$ is the $d$-dimensional forward component and $\mathbf{Y} = (Y_t : 0 \leq t \leq T)$ the one-dimensional backward one (the extension of our results to multidimensional backward equations is straightforward). Here, $W$ is a $q$-dimensional Brownian motion defined on a filtered probability space $(\Omega, \mathcal{F}, \mathbb{P}, (\mathcal{F}_t)_{0 \leq t \leq T})$, where $(\mathcal{F}_t)_t$ is the augmented natural filtration of $W$. The driver $f(\cdot, \cdot, \cdot, \cdot)$ and the terminal condition $\Phi(\cdot)$ are, respectively, a deterministic function and a deterministic functional of the process $\mathbf{S}$. The assumptions (H1)–(H3) below ensure the existence and the uniqueness of a solution $(\mathbf{S}, \mathbf{Y}, \mathbf{Z})$ to such equation (1)–(2).

Received June 2004; revised January 2005.
[1]Supported by Association Nationale de la Recherche Technique, École Polytechnique and Électricité de France.

*AMS 2000 subject classifications.* 60H10, 60H10, 65C30.

*Key words and phrases.* Backward stochastic differential equations, regression on function bases, Monte Carlo methods.

*Applications of BSDEs.* Such equations, first studied by Pardoux and Peng [26] in a general form, are important tools in mathematical finance. We mention some applications and refer the reader to [10, 12] for numerous references. In a complete market, for the usual valuation of a contingent claim with payoff $\Phi(\mathbf{S})$, $Y$ is the value of the replicating portfolio and $Z$ is related to the hedging strategy. In that case, the driver $f$ is linear w.r.t. $Y$ and $Z$. Some market imperfections can also be incorporated, such as higher interest rate for borrowing [4]: then, the driver is only Lipschitz continuous w.r.t. $Y$ and $Z$. Related numerical experiments are developed in Section 6. In incomplete markets, the Föllmer–Schweizer strategy [14] is given by the solution of a BSDE. When trading constraints on some assets are imposed, the super-replication price [13] is obtained as the limit of nonlinear BSDEs. Connections with recursive utilities of Duffie and Epstein [11] are also available. Peng has introduced the notion of $g$-expectation (here $g$ is the driver) as a nonlinear pricing rule [28]. Recently he has shown [27] the deep connection between BSDEs and dynamic risk measures, proving that any dynamic risk measure $(\mathcal{E}_t)_{0 \leq t \leq T}$ (satisfying some axiomatic conditions) is necessarily associated to a BSDE $(Y_t)_{0 \leq t \leq T}$ (the converse being known for years). The least we can say is that BSDEs are now inevitable tools in mathematical finance. Another indirect application may concern variance reduction techniques for the Monte Carlo computations of expectations, say $\mathbb{E}(\Phi)$ taking $f \equiv 0$. Indeed, $\int_0^T Z_s \, dW_s$ is the so-called martingale control variate (see [24], for instance). Finally, for applications to semi-linear PDEs, we refer to [25], among others.

The mathematical analysis of BSDE is now well understood (see [23] for recent references) and its numerical resolution has made recent progresses. However, even if several numerical methods have been proposed, they suffer of a high complexity in terms of computational time or are very costly in terms of computer memory. Thus, their uses in practice on real problems are difficult. Hence, it is still topical to devise more efficient algorithms. This article contributes in this direction by developing a simple approach, based on Monte Carlo regression on function bases. It is in the vein of the general regression approach of Bouchard and Touzi [6], but here it is actually much simpler because only one set of paths is used to evaluate all the regression operators. Consequently, the numerical implementation is easier and more efficient. In addition, we provide a full mathematical analysis of the influence of the parameters of the method.

*Numerical methods for BSDEs.* In the past decade, there have been several attempts to provide approximation schemes for BSDEs. First, Ma, Protter and Yong [22] propose the *four step scheme* to solve general FBSDEs, which requires the numerical resolution of a quasilinear parabolic PDE. In [2], Bally presents a time discretization scheme based on a Poisson net: this trick avoids him using the unknown regularity of $Z$ and enables him to derive a rate of convergence w.r.t. the intensity of the Poisson process. However, extra computations of very high-dimensional integrals are needed and this is not handled in [2]. In a recent work [29], Zhang proves some $\mathbf{L}_2$-regularity on $Z$, which allows the use of a regular deterministic time mesh. Under an assumption of *constructible functionals* for $\Phi$ (which essentially means that the system can be made Markovian, by adding $d'$ extra state variables), its approximation scheme is less consuming in terms of high-dimensional integrals. If for each of the $d+d'$ state variables, one uses $M$ points to compute the integrals, the complexity is about $M^{d+d'}$ per time step, for a global error of order $M^{-1}$ say (actually, an analysis of the global accuracy is not provided in [29]). This approach is somewhat related to the quantization method of Bally and Pagès [3], which is an optimal space discretization of the underlying dynamic programming equation (see also the former work by Chevance [8], where the driver does not depend on $Z$). We should also mention the works by Ma, Protter, San Martin and Soledad [21] and Briand, Delyon and Mémin [7], where the Brownian motion is replaced by a scaled random walk. Weak convergence results are given, without rates of approximation. The complexity becomes very large in multidimensional problems, like for finite differences schemes for PDEs. Recently, in the case of path-independent terminal conditions $\Phi(\mathbf{S}) = \phi(S_T)$, Bouchard and Touzi [6] propose a Monte Carlo approach which may be more suitable for high-dimensional problems. They follow the approach by Zhang [29] by approximating (1)–(2) by a discrete time FBSDE with $N$ time steps [see (5)–(6) below], with an $\mathbf{L}_2$-error of order $N^{-1/2}$. Instead of computing the conditional expectations which appear at each discretization time by discretizing the space of each state variable, the authors use a general regression operator, which can be derived, for instance, from kernel estimators or from the Malliavin calculus integration by parts formulas. The regression operator at a discretization time is assumed to be built independently of the underlying process, and independently of the regression operators at the other times. For the Malliavin calculus approach, for example, this means that one needs to simulate at each discrete time, $M$ copies of the approximation of (1), which is very costly. The algorithm that we propose in this paper requires only one set of paths to approximate all the regression operators at each discretization time at once. Since the regression operators are now correlated, the mathematical analysis is much more involved.

The regression operator we use in the sequel results from the $\mathbf{L}_2$-projection on a finite basis of functions, which leads in practice to solve a standard least squares problem. This approach is not new in numerical methods for financial engineering, since it has been developed by Longstaff and Schwartz [20] for the pricing of Bermuda options. See also [5] for the option pricing using simulations under the objective probability.

*Organization of the paper.* In Section 2 we set the framework of our study, define some notation used throughout the paper and describe our algorithm based on the approximation of conditional expectations by a projection on a finite basis of functions. We also provide some remarks related to models in finance.

The next three sections are devoted to analyzing the influence of the parameters of this scheme on the evaluation of $Y$ and $Z$. Note that approximation results on $Z$ were not previously considered in [6]. In Section 3 we provide an estimation of the time discretization error: this essentially follows from the results by Zhang [29]. Then, the impact of the function bases and the number of simulated paths is separately discussed in Section 4 and in Section 5, which is the major contribution of our work. Since this least squares approach is also popular to price Bermuda options [20], it is crucial to accurately estimate the propagation of errors in this type of numerical method, that is, to ensure that it is not explosive when the exercise frequency shrinks to 0. $\mathbf{L}_2$-estimates and a central limit theorem (see also [9] for Bermuda options) are proved.

In Section 6 explicit choices of function bases are given, together with numerical examples relative to the pricing of vanilla options and Asian options with differential interest rates.

## 2. Assumptions, notation and the numerical scheme.

2.1. *Standing assumptions.* Throughout the paper we assume that the following hypotheses are fulfilled:

(H1) The functions $(t,x) \mapsto b(t,x)$ and $(t,x) \mapsto \sigma(t,x)$ are uniformly Lipschitz continuous w.r.t. $(t,x) \in [0,T] \times \mathbb{R}^d$.
(H2) The driver $f$ satisfies the following continuity estimate:

$$|f(t_2, x_2, y_2, z_2) - f(t_1, x_1, y_1, z_1)|$$
$$\leq C_f(|t_2 - t_1|^{1/2} + |x_2 - x_1| + |y_2 - y_1| + |z_2 - z_1|)$$

for any $(t_1, x_1, y_1, z_1), (t_2, x_2, y_2, z_2) \in [0,T] \times \mathbb{R}^d \times \mathbb{R} \times \mathbb{R}^q$.
(H3) The terminal condition $\Phi$ satisfies the *functional Lipschitz condition*, that is, for any continuous functions $\mathbf{s}^1$ and $\mathbf{s}^2$, one has

$$|\Phi(\mathbf{s}^1) - \Phi(\mathbf{s}^2)| \leq C \sup_{t \in [0,T]} |s_t^1 - s_t^2|.$$

These assumptions (H1)–(H3) are sufficient to ensure the existence and uniqueness of a triplet $(\mathbf{S}, \mathbf{Y}, \mathbf{Z})$ solution to (1)–(2) (see [23] and references therein). In addition, the assumption (H3) allows a large class of terminal conditions (see examples in Section 2.4).

To approximate the forward component (1), we use a standard Euler scheme with time step $h$ (say smaller than 1), associated to equidistant discretization times $(t_k = kh = kT/N)_{0 \leq k \leq N}$. This approximation is defined by $S_0^N = S_0$ and

(3) $\qquad S_{t_{k+1}}^N = S_{t_k}^N + b(t_k, S_{t_k}^N)h + \sigma(t_k, S_{t_k}^N)(W_{t_{k+1}} - W_{t_k}).$

The terminal condition $\Phi(\mathbf{S})$ is approximated by $\Phi^N(P_{t_N}^N)$, where $\Phi^N$ is a deterministic function and $(P_{t_k}^N)_{0 \leq k \leq N}$ is a Markov chain, whose first components are given by those of $(S_{t_k}^N)_{0 \leq k \leq N}$. In other words, we eventually add extra state variables to make Markovian the implicit dynamics of the terminal condition. We also assume that $P_{t_k}^N$ is $\mathcal{F}_{t_k}$-measurable and that $\mathbb{E}[\Phi^N(P_{t_N}^N)]^2 < \infty$. Of course, this approximation strongly depends on the terminal condition type and its impact is measured by the error $\mathbb{E}|\Phi(\mathbf{S}) - \Phi^N(P_{t_N}^N)|^2$ (see Theorem 1 later). Examples of function $\Phi^N$ are given in Section 2.4.

Another hypothesis is required to prove that a certain discrete time BSDE $(Y_{t_k}^N)_k$ can be represented as a Lipschitz continuous function $y^N(t_k, \cdot)$ of $P_{t_k}^N$ (see Proposition 3 later). This property is mainly used in Section 6 on numerical experiments to derive relevant regression approximations.

(H4) The function $\Phi^N(\cdot)$ is Lipschitz continuous (uniformly in $N$) and $\sup_N |\Phi^N(\mathbf{0})| < \infty$. In addition, $\mathbb{E}|P_{t_N}^{N,k_0,x} - P_{t_N}^{N,k_0,x'}|^2 + \mathbb{E}|P_{t_{k_0+1}}^{N,k_0,x} - P_{t_{k_0+1}}^{N,k_0,x'}|^2 \leq C|x - x'|^2$ uniformly in $k_0$ and $N$.

Here, $(P_{t_k}^{N,k_0,x})_k$ stands for the Markov chain $(P_{t_k}^N)_k$ starting at $P_{t_{k_0}}^N = x$. Moreover, since we deal with the flow properties of $(P_{t_k}^N)_k$, we use the standard representation of this Markov chain as a random iterative sequence of the form $P_{t_k}^N = F_k^N(U_k, P_{t_{k-1}}^N)$, where $(F_k^N)_k$ are measurable functions and $(U_k)_k$ are i.i.d. random variables.

### 2.2. Notation.

PROJECTION ON FUNCTION BASES.

- The $\mathbf{L}_2(\Omega, \mathbb{P})$ projection of the random variable $U$ on a finite family $\phi = [\phi_1, \ldots, \phi_n]^*$ (considered as a random column vector) is denoted by $\mathcal{P}_\phi(U)$. We set $\mathcal{R}_\phi(U) = U - \mathcal{P}_\phi(U)$ for the projection error.
- At each time $t_k$, to approximate, respectively, $Y_{t_k}$ and $Z_{l,t_k}$ ($Z_{l,t_k}$ is the $l$th component of $Z_{t_k}$, $1 \leq l \leq q$), we will use, respectively, finite-dimensional function bases $p_{0,k}(P_{t_k}^N)$ and $p_{l,k}(P_{t_k}^N)$ ($1 \leq l \leq q$), which may be also written $p_{0,k}$ and $p_{l,k}$ ($1 \leq l \leq q$) to simplify. In the following, for convenience, both $(p_{l,k}(\cdot))$ and $(p_{l,k}(P_{t_k}^N))$ are indifferently called *function basis*. Explicit examples are given in Section 6. The projection coefficients will be denoted $\alpha_{0,k}, \alpha_{1,k}, \ldots, \alpha_{q,k}$ (viewed as column vectors). We assume that $\mathbb{E}|p_{l,k}|^2 < \infty$ ($0 \leq l \leq q$) and w.l.o.g. that $\mathbb{E}(p_{l,k}p_{l,k}^*)$ is invertible, which ensures the uniqueness of the coefficients of the projection $\mathcal{P}_{p_{l,k}}$ ($0 \leq l \leq q$).

- To simplify, we write $f_k(\alpha_{0,k}, \ldots, \alpha_{q,k})$ or $f_k(\alpha_k)$ for $f(t_k, S_{t_k}^N, \alpha_{0,k} \cdot p_{0,k}, \ldots, \alpha_{q,k} \cdot p_{q,k})$ [$S_{t_k}^N$ is the Euler approximation of $S_{t_k}$, see (3)].
- For convenience, we write $\mathbb{E}_k(\cdot) = \mathbb{E}(\cdot|\mathcal{F}_{t_k})$. We put $\Delta W_k = W_{t_{k+1}} - W_{t_k}$ (and $\Delta W_{l,k}$ component-wise) and define $v_k$ the (column) vector given by $[v_k]^* = (p_{0,k}^*, p_{1,k}^* \frac{\Delta W_{1,k}}{\sqrt{h}}, \ldots, p_{q,k}^* \frac{\Delta W_{q,k}}{\sqrt{h}})$.
- For a vector $x$, $|x|$ stands, as usual, for its Euclidean norm. The relative dimension is still implicit. For an integer $M$ and $x \in \mathbb{R}^M$, we put $|x|_M^2 = \frac{1}{M} \sum_{m=1}^M |x_m|^2$. For a set of projection coefficients $\alpha = (\alpha_0, \ldots, \alpha_q)$, we set $|\alpha| = \max_{0 \leq l \leq q} |\alpha_l|$ (the dimensions of the $\alpha_l$ may be different). For the set of basis functions at a fixed time $t_k$, $|p_k|$ is defined analogously.
- For a real symmetric matrix $A$, $\|A\|$ and $\|A\|_F$ are, respectively, the maximum of the absolute value of its eigenvalues and its Frobenius norm (defined by $\|A\|_F^2 = \sum_{i,j} a_{i,j}^2$).

We refer to Section 6 for explicit choices of function bases, but to fix ideas, a possible choice could be to define, for each time $t_k$, grids $(x_{l,k}^i : 1 \leq i \leq n)_{0 \leq l \leq q}$ and define $p_{l,k}(\cdot)$ as the basis of indicator functions of the open Voronoi partition [17] associated to $(x_{l,k}^i : 1 \leq i \leq n)$, that is, $p_{l,k}(\cdot) = (\mathbf{1}_{C_{l,k}^i}(\cdot))_{1 \leq i \leq n}$, where $C_{l,k}^i = \{x : |x - x_{l,k}^i| < |x - x_{l,k}^j|, \forall j \neq i\}$.

SIMULATIONS. In the following, $M$ independent simulations of $(P_{t_k}^N)_{0 \leq k \leq N}$, $(\Delta W_k)_{0 \leq k \leq N-1}$ will be used. We denote them $((P_{t_k}^{N,m})_{0 \leq k \leq N})_{1 \leq m \leq M}$, $((\Delta W_k^m)_{0 \leq k \leq N-1})_{1 \leq m \leq M}$:

- The values of basis functions along these simulations are denoted $(p_{l,k}^m = p_l(P_{t_k}^{N,m}))_{0 \leq l \leq q, 0 \leq k \leq N-1, 1 \leq m \leq M}$.
- Analogously to $f_k(\alpha_{0,k}, \ldots, \alpha_{q,k})$ or $f_k(\alpha_k)$, we denote $f_k^m(\alpha_{0,k}, \ldots, \alpha_{q,k})$ or $f_k^m(\alpha_k)$ for $f(t_k, S_{t_k}^{N,m}, \alpha_{0,k} \cdot p_{0,k}^m, \ldots, \alpha_{q,k} \cdot p_{q,k}^m)$.

We define the following:

- the (column) vector $v_k^m$ by $[v_k^m]^* = (p_{0,k}^{m*}, p_{1,k}^{m*} \frac{\Delta W_{1,k}^m}{\sqrt{h}}, \ldots, p_{q,k}^{m*} \frac{\Delta W_{q,k}^m}{\sqrt{h}})$;
- the matrix $V_k^M = \frac{1}{M} \sum_{m=1}^M v_k^m [v_k^m]^*$;
- the matrix $P_{l,k}^M = \frac{1}{M} \sum_{m=1}^M p_{l,k}^m [p_{l,k}^m]^*$ ($0 \leq l \leq q$).

TRUNCATIONS. To ensure the stability of the algorithm, we use threshold techniques, which are based on the following notation:

- In Proposition 2 below, based on BSDEs' a priori estimates, we explicitly build some $\mathbb{R}$-valued functions $(\rho_{l,k}^N)_{0 \leq l \leq q, 0 \leq k \leq N-1}$ bounded from below by 1. We set $\rho_k^N(P_{t_k}^N) = [\rho_{0,k}^N(P_{t_k}^N), \ldots, \rho_{q,k}^N(P_{t_k}^N)]^*$.

- Associated to these estimates, we define (random) truncation functions $\hat{\rho}_{l,k}^N(x) = \rho_{l,k}^N(P_{t_k}^N)\xi(x/\rho_{l,k}^N(P_{t_k}^N))$ and $\hat{\rho}_{l,k}^{N,m}(x) = \rho_{l,k}^N(P_{t_k}^{N,m})\xi(x/\rho_{l,k}^N(P_{t_k}^{N,m}))$, where $\xi:\mathbb{R} \mapsto \mathbb{R}$ is a $C_b^2$-function, such that $\xi(x) = x$ for $|x| \leq 3/2$, $|\xi|_\infty \leq 2$ and $|\xi'|_\infty \leq 1$.

In the next computations, $C$ denotes a generic constant that may change from line to line. It is still uniform in the parameters of our scheme.

2.3. *The numerical scheme.* We are now in a position to define the simulation-based approximations of the BSDE (1)–(2). The statements of approximation results and their proofs are postponed to Sections 3, 4 and 5.

Our procedure combines a backward in time evaluation (from time $t_N = T$ to time $t_0 = 0$), a fixed point argument (using $i = 1, \ldots, I$ Picard iterations), least squares problems on $M$ simulated paths (using some function bases).

*Initialization.* The algorithm is initialized with $Y_{t_N}^{N,i,I,M} = \Phi^N(P_{t_N}^N)$ (independently of $i$ and $I$). Then, the solution $(Y_{t_k}, Z_{1,t_k}, \ldots, Z_{q,t_k})$ at a given time $t_k$ is represented via some projection coefficients $(\alpha_{l,k}^{i,I,M})_{0 \leq l \leq q}$ by

$$Y_{t_k}^{N,i,I,M} = \hat{\rho}_{0,k}^N(\alpha_{0,k}^{i,I,M} \cdot p_{0,k}), \qquad \sqrt{h} Z_{l,t_k}^{N,i,I,M} = \hat{\rho}_{l,k}^N(\sqrt{h}\alpha_{l,k}^{i,I,M} \cdot p_{l,k})$$

($\hat{\rho}_{0,k}^N$ and $\hat{\rho}_{l,k}^N$ are the truncations introduced before). We now detail how the coefficients are computed using independent realizations $((P_{t_k}^{N,m})_{0 \leq k \leq N})_{1 \leq m \leq M}$, $((\Delta W_k^m)_{0 \leq k \leq N-1})_{1 \leq m \leq M}$.

*Backward in time iteration at time $t_k < T$.* Assume that an approximation $Y_{t_{k+1}}^{N,I,I,M} := \hat{\rho}_{0,k+1}^N(\alpha_{0,k+1}^{I,I,M} \cdot p_{0,k+1})$ is built, and denote $Y_{t_{k+1}}^{N,I,I,M,m} = \hat{\rho}_{0,k+1}^{N,m}(\alpha_{0,k+1}^{I,I,M} \cdot p_{0,k+1}^m)$ its realization along the $m$th simulation.

→ For the initialization $i = 0$ of Picard iterations, set $Y_{t_k}^{N,0,I,M} = 0$ and $Z_{t_k}^{N,0,I,M} = 0$, that is, $\alpha_{l,k}^{0,I,M} = 0$ ($0 \leq l \leq q$).
→ For $i = 1, \ldots, I$, the coefficients $\alpha_k^{i,I,M} = (\alpha_{l,k}^{i,I,M})_{0 \leq l \leq q}$ are iteratively obtained as the arg min in $(\alpha_0, \ldots, \alpha_q)$ of the quantity

$$(4) \quad \frac{1}{M} \sum_{m=1}^{M} \left( Y_{t_{k+1}}^{N,I,I,M,m} - \alpha_0 \cdot p_{0,k}^m + h f_k^m(\alpha_k^{i-1,I,M}) - \sum_{l=1}^{q} \alpha_l \cdot p_{l,k}^m \Delta W_{l,k}^m \right)^2.$$

If the above least squares problem has multiple solutions (i.e., the empirical regression matrix is not invertible, which occurs with small probability when $M$ becomes large), we may choose, for instance, the (unique) solution of minimal norm. Actually, this choice is arbitrary and has no incidence on the further analysis.

The convergence parameters of this scheme are the time step $h$ ($h \to 0$), the function bases, the number of simulations $M$ ($M \to +\infty$). This is fully analyzed in the following sections, with three main steps: time discretization of the BSDE, projections on bases functions in $\mathbf{L}_2(\Omega, \mathbb{P})$, empirical projections using simulated paths. An estimate of the global error directly follows from the combination of Theorems 1, 2 and 3. We will also see that it is enough to have $I = 3$ Picard iterations (see Theorem 3).

The intuition behind the above sequence of least squares problems (4) is actually simple. It aims at mimicking what can be ideally done with an infinite number of simulations, Picard iterations and bases functions, that is,

$$(Y^N_{t_k}, Z^N_{t_k}) = \mathop{\arg\inf}_{(Y,Z) \in \mathbf{L}_2(\mathcal{F}_{t_k})} \mathbb{E}(Y^N_{t_{k+1}} - Y + hf(t_k, S^N_{t_k}, Y, Z) - Z\Delta W_k)^2,$$

where, as usual, $\mathbf{L}_2(\mathcal{F}_{t_k})$ stands for the square integrable and $\mathcal{F}_{t_k}$-measurable, possibly multidimensional, random variables. This ideal case is an appoximation of the BSDE (2) which writes

$$Y_{t_{k+1}} + \int_{t_k}^{t_{k+1}} f(s, S_s, Y_s, Z_s)\, ds = Y_{t_k} + \int_{t_k}^{t_{k+1}} Z_s\, dW_s$$

over the time interval $[t_k, t_{k+1}]$. $(Y^N_{t_k})_k$ will be interpreted as a discrete time BSDE (see Theorem 1).

2.4. *Remarks for models in finance.* Here, we give examples of drivers $f$ and terminal conditions $\Phi(\mathbf{S})$ in the case of option pricing with different interest rates [4]: $R$ for borrowing and $r$ for lending with $R \geq r$. Assume for simplicity that there is only one underlying risky asset ($d = 1$) whose dynamics is given by the Black–Scholes model with drift $\mu$ and volatility $\sigma$ ($q = 1$): $dS_t = S_t(\mu\, dt + \sigma\, dW_t)$.

- *Driver*: If we set $f(t, x, y, z) = -\{yr + z\theta - (y - \frac{z}{\sigma})^-(R - r)\}$, where $\theta = \frac{\mu - r}{\sigma}$, $Y_t$ is the value at time $t$ of the self-financing portfolio replicating the payoff $\Phi(\mathbf{S})$ [12]. In the case of equal interest rates $R = r$, the driver is linear and we obtain the usual risk-neutral valuation rule.
- *Terminal conditions*: A large class of exotic payoffs satisfies the functional Lipschitz condition (H3).
  - Vanilla payoff: $\Phi(\mathbf{S}) = \phi(S_T)$. Set $P^N_{t_k} = S^N_{t_k}$ and $\Phi^N(P^N_{t_N}) = \phi(P^N_{t_N})$. Under (H3), it gives $\mathbb{E}|\Phi^N(P^N_{t_N}) - \Phi(\mathbf{S})|^2 \leq Ch$.
  - Asian payoff: $\Phi(\mathbf{S}) = \phi(S_T, \int_0^T S_t\, dt)$. Set $P^N_{t_k} = (S^N_{t_k}, h\sum_{i=0}^{k-1} S^N_{t_i})$ and $\Phi^N(P^N_{t_N}) = \phi(P^N_{t_N})$. For usual functions $\phi$, the $\mathbf{L}_2$-error is of order $1/2$ w.r.t. $h$. More accurate approximations of the average of $\mathbf{S}$ could be incorporated [18].

- Lookback payoff: $\Phi(\mathbf{S}) = \phi(S_T, \min_{t\in[0,T]} S_t, \max_{t\in[0,T]} S_t)$. Set $\Phi^N(P^N_{t_N}) = \phi(P^N_{t_N})$ with $P^N_{t_k} = (S^N_{t_k}, \min_{i\leq k} S^N_{t_i}, \max_{i\leq k} S^N_{t_i})$. In general, this induces an $\mathbf{L}_2$-error of magnitude $\sqrt{h\log(1/h)}$ [29]. The rate $\sqrt{h}$ can be achieved by considering the exact extrema of the continuous Euler scheme [1].

Note also that (H4) is satisfied on these payoffs.

We also mention that the price process $(S_t)_t$ is usually positive coordinate-wise, but its Euler scheme [defined in (3)] does not enjoy this feature. This may be an undesirable property, which can be avoided by considering the Euler scheme on the log-price. With this modification, the analysis below is unchanged and we refer to [15] for details.

**3. Approximation results: step 1.** We first consider a time approximation of equations (1) and (2). The forward component is approximated using the Euler scheme (3) and the backward component (2) is evaluated in a backward manner. First, we set $Y^N_{t_N} = \Phi^N(P^N_{t_N})$. Then, $(Y^N_{t_k}, Z^N_{t_k})_{0\leq k\leq N-1}$ are defined by

$$Z^N_{l,t_k} = \frac{1}{h}\mathbb{E}_k(Y^N_{t_{k+1}}\Delta W_{l,k}), \tag{5}$$

$$Y^N_{t_k} = \mathbb{E}_k(Y^N_{t_{k+1}}) + hf(t_k, S^N_{t_k}, Y^N_{t_k}, Z^N_{t_k}). \tag{6}$$

Using, in particular, the inequality $|Z^N_{l,t_k}| \leq \frac{1}{\sqrt{h}}\sqrt{\mathbb{E}_k(Y^N_{t_{k+1}})^2}$, it is easy to see by a recursive argument that $Y^N_{t_k}$ and $Z^N_{t_k}$ belong to $\mathbf{L}_2(\mathcal{F}_{t_k})$. It is also equivalent to assert that they minimize the quantity

$$\mathbb{E}(Y^N_{t_{k+1}} - Y + hf(t_k, S^N_{t_k}, Y, Z) - Z\Delta W_k)^2 \tag{7}$$

over $\mathbf{L}_2(\mathcal{F}_{t_k})$ random variables $(Y, Z)$. Note that $Y^N_{t_k}$ is well defined in (6), because the mapping $Y \mapsto \mathbb{E}_k(Y^N_{t_{k+1}}) + hf(t_k, S^N_{t_k}, Y, Z^N_{t_k})$ is a contraction in $\mathbf{L}_2(\mathcal{F}_{t_k})$, for $h$ small enough. The following result provides an estimate of the error induced by this first step.

THEOREM 1. *Assume* (H1)–(H3). *For $h$ small enough, we have*

$$\max_{0\leq k\leq N} \mathbb{E}|Y_{t_k} - Y^N_{t_k}|^2 + \sum_{k=0}^{N-1}\int_{t_k}^{t_{k+1}}\mathbb{E}|Z_t - Z^N_{t_k}|^2\,dt$$
$$\leq C((1+|S_0|^2)h + \mathbb{E}|\Phi(\mathbf{S}) - \Phi^N(P^N_{t_N})|^2).$$

PROOF. From [29], we know that the key point is the $\mathbf{L}^2$-regularity of $Z$. Here, under (H1)–(H3), $Z$ is càdlàg (see Remark 2.6.ii in [29]). Thus, The-orem 3.1 in [29] states that

$$\sum_{k=0}^{N-1} \mathbb{E} \int_{t_k}^{t_{k+1}} |Z_t - Z_{t_k}|^2 \, dt \le C(1 + |S_0|^2)h.$$

With this estimate, the proof of Theorem 1 is standard (see, e.g., the proof of Theorem 5.3 in [29]) and we omit details. $\square$

Owing to the Markov chain $(P_{t_k}^N)_{0 \le k \le N}$, the independent increments $(\Delta W_k)_{0 \le k \le N-1}$ and (5)–(6), we easily get the following result.

PROPOSITION 1. *Assume* (H1)–(H3). *For $h$ small enough, we have*

$$(8) \quad Y_{t_k}^N = y_k^N(P_{t_k}^N), \qquad Z_{l,t_k}^N = z_{l,k}^N(P_{t_k}^N) \qquad \text{for } 0 \le k \le N \text{ and } 1 \le l \le q,$$

*where* $(y_k^N(\cdot))_k$ *and* $(z_{l,k}^N(\cdot))_{k,l}$ *are measurable functions.*

It will be established in Section 6 that they are Lipschitz continuous under the extra assumption (H4).

**4. Approximation results: step 2.** Here, the conditional expectations which appear in the definitions (5)–(6) of $Y_{t_k}^N$ and $Z_{l,t_k}^N$ ($1 \le l \le q$) are replaced by a $\mathbf{L}_2(\Omega, \mathbb{P})$ projection on the function bases $p_{0,k}$ and $p_{l,k}$ ($1 \le l \le q$). A numerical difficulty still remains in the approximation of $Y_{t_k}^N$ in (6), which is usually obtained as a fixed point. To circumvent this problem, we propose a solution combining the projection on the function basis and $I$ Picard iterations. The integer $I$ is a fixed parameter of our scheme (the analysis below shows that the value $I = 3$ is relevant).

DEFINITION 1. We denote by $Y_{t_k}^{N,i,I}$ the approximation of $Y_{t_k}^N$, where $i$ Picard iterations with projections have been performed at time $t_k$ and $I$ Picard iterations with projections at any time after $t_k$. Analogous notation stands for $Z_{l,t_k}^{N,i,I}$. We associate to $Y_{t_k}^{N,i,I}$ and $Z_{l,t_k}^{N,i,I}$ their respective projection coefficients $\alpha_{0,k}^{i,I}$ and $\alpha_{l,k}^{i,I}$, on the function bases $p_{0,k}$ and $p_{l,k}$ ($1 \le l \le q$).

We now turn to a precise definition of the above quantities. We set $Y_{t_N}^{N,i,I} = \Phi^N(P_{t_N}^N)$, independently of $i$ and $I$. Assume that $Y_{t_{k+1}}^{N,I,I}$ is obtained and let us define $Y_{t_k}^{N,i,I}, Z_{l,t_k}^{N,i,I}$ for $i = 0, \dots, I$. We begin with $Y_{t_k}^{N,0,I} = 0$ and $Z_{t_k}^{N,0,I} = 0$, corresponding to $\alpha_{l,k}^{0,I} = 0$ ($0 \le l \le q$). By analogy with (7), we set $\alpha_k^{i,I} = (\alpha_{l,k}^{i,I})_{0 \le l \le q}$ as the arg min in $(\alpha_0, \dots, \alpha_q)$ of the quantity

$$(9) \quad \mathbb{E}\left(Y_{t_{k+1}}^{N,I,I} - \alpha_0 \cdot p_{0,k} + h f_k(\alpha_k^{i-1,I}) - \sum_{l=1}^{q} \alpha_l \cdot p_{l,k} \Delta W_{l,k}\right)^2.$$

Iterating with $i = 1, \ldots, I$, at the end we get $(\alpha_{l,k}^{I,I})_{0 \leq l \leq q}$, thus, $Y_{t_k}^{N,I,I} = \alpha_{0,k}^{I,I} \cdot p_{0,k}$ and $Z_{l,t_k}^{N,I,I} = \alpha_{l,k}^{I,I} \cdot p_{l,k}$ ($1 \leq l \leq q$). The least squares problem (9) can be formulated in different ways but this one is more convenient to get an intuition on (4). The error induced by this second step is analyzed by the following result.

THEOREM 2. *Assume* (H1)–(H3). *For $h$ small enough, we have*

$$\max_{0 \leq k \leq N} \mathbb{E}|Y_{t_k}^{N,I,I} - Y_{t_k}^N|^2 + h \sum_{k=0}^{N-1} \mathbb{E}|Z_{t_k}^{N,I,I} - Z_{t_k}^N|^2$$

$$\leq Ch^{2I-2}[1 + |S_0|^2 + \mathbb{E}|\Phi^N(P_{t_N}^N)|^2]$$

$$+ C \sum_{k=0}^{N-1} \mathbb{E}|\mathcal{R}_{p_{0,k}}(Y_{t_k}^N)|^2 + Ch \sum_{k=0}^{N-1} \sum_{l=1}^{q} \mathbb{E}|\mathcal{R}_{p_{l,k}}(Z_{l,t_k}^N)|^2.$$

The above result shows how projection errors cumulate along the backward iteration. The key point is to note that they only sum up, with a factor $C$ which does not explode as $N \to \infty$. These estimates improve those of Theorem 4.1 in [6] for two reasons. First, error estimates on $Z^N$ are provided here. Second, in the cited theorem, the error is analyzed in terms of $\mathbb{E}|\mathcal{R}_{p_{0,k}}(Y_{t_k}^{N,I,I})|^2$ and $\mathbb{E}|\mathcal{R}_{p_{l,k}}(Z_{l,t_k}^{N,I,I})|^2$ say: hence, the influence of function bases is still questionable, since it is hidden in the projection residuals $\mathcal{R}_{p_k}$ and also in the random variables $Y_{t_k}^{N,I,I}$ and $Z_{l,t_k}^{N,I,I}$. Our estimates are relevant to directly analyze the influence of function bases (see Section 6 for explicit computations). This feature is crucial in our opinion. Regarding the influence of $I$, it is enough here to have $I = 2$ to get an error of the same order as in Theorem 1. At the third step, $I = 3$ is needed.

PROOF OF THEOREM 2. For convenience, we denote $\mathcal{A}^N(S_0) = 1 + |S_0|^2 + \mathbb{E}|\Phi^N(P_{t_N}^N)|^2$. In the following computations, we repeatedly use three standard inequalities:

1. The contraction property of the $\mathbf{L}_2$-projection operator: for any random variable $X \in \mathbf{L}_2$, we have $\mathbb{E}|\mathcal{P}_{p_{l,k}}(X)|^2 \leq \mathbb{E}|X|^2$.
2. The Young inequality: $\forall \gamma > 0$, $\forall (a,b) \in \mathbb{R}^2$, $(a+b)^2 \leq (1+\gamma h)a^2 + (1+\frac{1}{\gamma h})b^2$.
3. The discrete Gronwall lemma: for any nonnegative sequences $(a_k)_{0 \leq k \leq N}$, $(b_k)_{0 \leq k \leq N}$ and $(c_k)_{0 \leq k \leq N}$ satisfying $a_{k-1} + c_{k-1} \leq (1+\gamma h)a_k + b_{k-1}$ (with $\gamma > 0$), we have $a_k + \sum_{i=k}^{N-1} c_i \leq e^{\gamma(T-t_k)}[a_N + \sum_{i=k}^{N-1} b_i]$. Most of the time, it will be used with $c_i = 0$.

Because $\Delta W_k$ is centered and independent of $(p_{l,k})_{0\leq l\leq q}$, it is straightforward to see that the solution of the least squares problem (9) is given, for $i \geq 1$, by

$$Z_{l,t_k}^{N,i,I} = \frac{1}{h}\mathcal{P}_{p_{l,k}}(Y_{t_{k+1}}^{N,I,I}\Delta W_{l,k}), \tag{10}$$

$$Y_{t_k}^{N,i,I} = \mathcal{P}_{p_{0,k}}(Y_{t_{k+1}}^{N,I,I} + hf(t_k, S_{t_k}^N, Y_{t_k}^{N,i-1,I}, Z_{t_k}^{N,i-1,I})). \tag{11}$$

The proof of Theorem 2 may be divided in several steps.

*Step* 1: *a* (*tight*) *preliminary upper bound for* $\mathbb{E}|Z_{l,t_k}^{N,i,I}|^2$. First note that $Z_{l,t_k}^{N,i,I}$ is constant for $i \geq 1$. Moreover, the Cauchy–Schwarz inequality yields $|\mathbb{E}_k(Y_{t_{k+1}}^{N,I,I}\Delta W_{l,k})|^2 = |\mathbb{E}_k([Y_{t_{k+1}}^{N,I,I} - \mathbb{E}_k(Y_{t_{k+1}}^{N,I,I})]\Delta W_{l,k})|^2 \leq h(\mathbb{E}_k[Y_{t_{k+1}}^{N,I,I}]^2 - [\mathbb{E}_k(Y_{t_{k+1}}^{N,I,I})]^2)$. Since $(p_{l,k})_l$ is $\mathcal{F}_{t_k}$-measurable and owing to the contraction of the projection operator, it follows that

$$\mathbb{E}|Z_{l,t_k}^{N,i,I}|^2 = \frac{1}{h^2}\mathbb{E}[\mathcal{P}_{p_{l,k}}(\mathbb{E}_k[Y_{t_{k+1}}^{N,I,I}\Delta W_{l,k}])]^2 \leq \frac{1}{h^2}\mathbb{E}(\mathbb{E}_k[Y_{t_{k+1}}^{N,I,I}\Delta W_{l,k}])^2$$
$$\leq \frac{1}{h}(\mathbb{E}[Y_{t_{k+1}}^{N,I,I}]^2 - \mathbb{E}[\mathbb{E}_k(Y_{t_{k+1}}^{N,I,I})]^2). \tag{12}$$

As it may be seen in the computations below, the term $\mathbb{E}[\mathbb{E}_k(Y_{t_{k+1}}^{N,I,I})]^2$ in (12) plays a crucial role to make further estimates not explosive w.r.t. $h$.

*Step* 2: $\mathbf{L}_2$ *bounds for* $Y_{t_k}^{N,i,I}$ *and* $\sqrt{h}Z_{l,t_k}^{N,i,I}$. Actually, it is an easy exercise to check that the random variables $Y_{t_k}^{N,i,I}$ and $\sqrt{h}Z_{l,t_k}^{N,i,I}$ are square integrable. We aim at proving that uniform $\mathbf{L}_2$ bounds w.r.t. $i, I, k$ are available. Denote $\chi_k^{N,I} : Y \in \mathbf{L}_2(\mathcal{F}_{t_k}) \mapsto \mathcal{P}_{p_{0,k}}(Y_{t_{k+1}}^{N,I,I} + hf(t_k, S_{t_k}^N, Y, Z_{t_k}^{N,i-1,I})) \in \mathbf{L}_2(\mathcal{F}_{t_k})$. Clearly, $\mathbb{E}|\chi_k^{N,I}(Y_2) - \chi_k^{N,I}(Y_1)|^2 \leq (C_f h)^2 \mathbb{E}|Y_2 - Y_1|^2$, where $C_f$ is the Lipschitz constant of $f$. Consequently, for $h$ small enough, the application $\chi_k^{N,I}$ is contracting and has a unique fixed point $Y_{t_k}^{N,\infty,I} \in \mathbf{L}_2(\mathcal{F}_{t_k})$ (remind that $Z_{l,t_k}^{N,i,I}$ does not depend on $i \geq 1$). One has

$$Y_{t_k}^{N,\infty,I} = \mathcal{P}_{p_{0,k}}(Y_{t_{k+1}}^{N,I,I} + hf(t_k, S_{t_k}^N, Y_{t_k}^{N,\infty,I}, Z_{t_k}^{N,I,I})), \tag{13}$$

$$\mathbb{E}|Y_{t_k}^{N,\infty,I} - Y_{t_k}^{N,i,I}|^2 \leq (C_f h)^{2i}\mathbb{E}|Y_{t_k}^{N,\infty,I}|^2 \tag{14}$$

since $Y_{t_k}^{N,0,I} = 0$. Thus, Young's inequality yields, for $i \geq 1$,

$$\mathbb{E}|Y_{t_k}^{N,i,I}|^2 \leq \left(1 + \frac{1}{h}\right)\mathbb{E}|Y_{t_k}^{N,\infty,I} - Y_{t_k}^{N,i,I}|^2 + (1+h)\mathbb{E}|Y_{t_k}^{N,\infty,I}|^2$$
$$\leq (1+Ch)\mathbb{E}|Y_{t_k}^{N,\infty,I}|^2. \tag{15}$$

The above inequality is also true for $i=0$ because $Y_{t_k}^{N,0,I}=0$. We now estimate $\mathbb{E}|Y_{t_k}^{N,\infty,I}|^2$ from the identity (13). Combining Young's inequality (with $\gamma$ to be chosen later), the identity $\mathcal{P}_{p_{0,k}}(Y_{t_{k+1}}^{N,I,I}) = \mathcal{P}_{p_{0,k}}(\mathbb{E}_k[Y_{t_{k+1}}^{N,I,I}])$, the contraction of $\mathcal{P}_{p_{0,k}}$ and the Lipschitz property of $f$, we get

$$\mathbb{E}|Y_{t_k}^{N,\infty,I}|^2 \leq (1+\gamma h)\mathbb{E}|\mathbb{E}_k[Y_{t_{k+1}}^{N,I,I}]|^2 \tag{16}$$
$$+ Ch\left(h + \frac{1}{\gamma}\right)[\mathbb{E}f_k^2(0,\ldots,0) + \mathbb{E}|Y_{t_k}^{N,\infty,I}|^2 + \mathbb{E}|Z_{t_k}^{N,I,I}|^2].$$

Bringing together terms $\mathbb{E}|Y_{t_k}^{N,\infty,I}|^2$, then using (12) and the easy upper bound $\mathbb{E}f_k^2(0,\ldots,0) \leq C(1+|S_0|^2)$, it readily follows that

$$\mathbb{E}|Y_{t_k}^{N,\infty,I}|^2 \leq \frac{(1+\gamma h)}{1-Ch(h+1/\gamma)}\mathbb{E}|\mathbb{E}_k[Y_{t_{k+1}}^{N,I,I}]|^2$$
$$\tag{17} + \frac{Ch(h+1/\gamma)}{1-Ch(h+1/\gamma)}[1+|S_0|^2]$$
$$+ \frac{C(h+1/\gamma)}{1-Ch(h+1/\gamma)}(\mathbb{E}|Y_{t_{k+1}}^{N,I,I}|^2 - \mathbb{E}|\mathbb{E}_k[Y_{t_{k+1}}^{N,I,I}]|^2),$$

provided that $h$ is small enough. Take $\gamma = C$ to get

$$\mathbb{E}|Y_{t_k}^{N,\infty,I}|^2 \leq Ch[1+|S_0|^2] + (1+Ch)\mathbb{E}|Y_{t_{k+1}}^{N,I,I}|^2 + Ch\mathbb{E}|\mathbb{E}_k[Y_{t_{k+1}}^{N,I,I}]|^2$$
$$\tag{18} \leq Ch[1+|S_0|^2] + (1+2Ch)\mathbb{E}|Y_{t_{k+1}}^{N,I,I}|^2$$

with a new constant $C$. Plugging this estimate into (15) with $i=I$, we get $\mathbb{E}|Y_{t_k}^{N,I,I}|^2 \leq Ch[1+|S_0|^2] + (1+Ch)\mathbb{E}|Y_{t_{k+1}}^{N,I,I}|^2$ and, thus, by Gronwall's lemma, $\sup_{0\leq k\leq N}\mathbb{E}|Y_{t_k}^{N,I,I}|^2 \leq C\mathcal{A}^N(S_0)$. This upper bound combined with (18), (15) and (12) finally provides the required uniform estimates for $\mathbb{E}|Y_{t_k}^{N,i,I}|^2$ and $\mathbb{E}|Z_{l,t_k}^{N,i,I}|^2$:

$$\sup_{I\geq 1}\sup_{i\geq 0}\sup_{0\leq k\leq N}(\mathbb{E}|Y_{t_k}^{N,i,I}|^2 + h\mathbb{E}|Z_{l,t_k}^{N,i,I}|^2) \leq C\mathcal{A}^N(S_0). \tag{19}$$

*Step* 3: *upper bounds for* $\eta_k^{N,I} = \mathbb{E}|Y_{t_k}^{N,I,I} - Y_{t_k}^N|^2$. Note that $\eta_N^{N,I}=0$. Our purpose is to prove the following relation for $0\leq k < N$:

$$\eta_k^{N,I} \leq (1+Ch)\eta_{k+1}^{N,I} + Ch^{2I-1}\mathcal{A}^N(S_0)$$
$$\tag{20} + C\mathbb{E}|\mathcal{R}_{p_{0,k}}(Y_{t_k}^N)|^2 + Ch\sum_{l=1}^{q}\mathbb{E}|\mathcal{R}_{p_{l,k}}(Z_{l,t_k}^N)|^2.$$

Note that the estimate on $\max_{0 \leq k \leq N} \mathbb{E}|Y_{t_k}^{N,I,I} - Y_{t_k}^N|^2$ given in Theorem 2 directly follows from the relation above. With the arguments used to derive (15) and using the estimate (19), we easily get

$$
\begin{aligned}
\eta_k^{N,I} &\leq Ch^{2I-1}\mathcal{A}^N(S_0) + (1+h)\mathbb{E}|Y_{t_k}^{N,\infty,I} - Y_{t_k}^N|^2 \\
&= Ch^{2I-1}\mathcal{A}^N(S_0) + (1+h)\mathbb{E}|\mathcal{R}_{p_{0,k}}(Y_{t_k}^N)|^2 \\
&\quad + (1+h)\mathbb{E}|Y_{t_k}^{N,\infty,I} - \mathcal{P}_{p_{0,k}}(Y_{t_k}^N)|^2,
\end{aligned}
\tag{21}
$$

where we used at the last equality the orthogonality property relative to $\mathcal{P}_{p_{0,k}}$:

$$
\mathbb{E}|Y_{t_k}^{N,\infty,I} - Y_{t_k}^N|^2 = \mathbb{E}|\mathcal{R}_{p_{0,k}}(Y_{t_k}^N)|^2 + \mathbb{E}|Y_{t_k}^{N,\infty,I} - \mathcal{P}_{p_{0,k}}(Y_{t_k}^N)|^2.
\tag{22}
$$

Furthermore, with the same techniques as for (12) and (16), we can prove

$$
\begin{aligned}
&\mathbb{E}|Z_{t_k}^{N,I,I} - Z_{t_k}^N|^2 \\
&= \sum_{l=1}^{q} \mathbb{E}|\mathcal{R}_{p_{l,k}}(Z_{l,t_k}^N)|^2 + \sum_{l=1}^{q} \mathbb{E}|Z_{l,t_k}^{N,I,I} - \mathcal{P}_{p_{l,k}}(Z_{l,t_k}^N)]|^2 \\
&\leq \sum_{l=1}^{q} \mathbb{E}|\mathcal{R}_{p_{l,k}}(Z_{l,t_k}^N)|^2 \\
&\quad + \frac{d}{h}(\mathbb{E}[Y_{t_{k+1}}^{N,I,I} - Y_{t_{k+1}}^N]^2 - \mathbb{E}[\mathbb{E}_k(Y_{t_{k+1}}^{N,I,I} - Y_{t_{k+1}}^N)]^2)
\end{aligned}
\tag{23}
$$

and

$$
\begin{aligned}
&\mathbb{E}|Y_{t_k}^{N,\infty,I} - \mathcal{P}_{p_{0,k}}(Y_{t_k}^N)|^2 \\
&\leq (1+\gamma h)\mathbb{E}|\mathbb{E}_k[Y_{t_{k+1}}^{N,I,I} - Y_{t_{k+1}}^N]|^2 \\
&\quad + Ch\left(h + \frac{1}{\gamma}\right)[\mathbb{E}|Y_{t_k}^{N,\infty,I} - Y_{t_k}^N|^2 + \mathbb{E}|Z_{t_k}^{N,I,I} - Z_{t_k}^N|^2].
\end{aligned}
\tag{24}
$$

Replacing the estimate (23) in (24), choosing $\gamma = Cd$ and using (22) directly leads to

$$
\begin{aligned}
&(1 - Ch)\mathbb{E}|Y_{t_k}^{N,\infty,I} - \mathcal{P}_{p_{0,k}}(Y_{t_k}^N)|^2 \\
&\leq (1 + Ch)\eta_{k+1}^{N,I} \\
&\quad + Ch\sum_{l=1}^{q}\mathbb{E}|\mathcal{R}_{p_{l,k}}(Z_{l,t_k}^N)|^2 + Ch\mathbb{E}|\mathcal{R}_{p_{0,k}}(Y_{t_k}^N)|^2.
\end{aligned}
\tag{25}
$$

Plugging this estimate into (21) completes the proof of (20). □

*Step* 4: *upper bounds for* $\zeta^N = h\sum_{k=0}^{N-1}\mathbb{E}|Z_{t_k}^{N,I,I} - Z_{t_k}^N|^2$. We aim at showing

$$\zeta^N \leq Ch^{2I-2}\mathcal{A}^N(S_0) + Ch\sum_{k=0}^{N-1}\sum_{l=1}^{q}\mathbb{E}|\mathcal{R}_{p_{l,k}}(Z_{l,t_k}^N)|^2$$

(26)
$$+ C\sum_{k=0}^{N-1}\mathbb{E}|\mathcal{R}_{p_{0,k}}(Y_{t_k}^N)|^2 + C\max_{0\leq k\leq N-1}\eta_k^{N,I}.$$

In view of (23), we have

$$\zeta^N \leq h\sum_{k=0}^{N-1}\sum_{l=1}^{q}\mathbb{E}|\mathcal{R}_{p_{l,k}}(Z_{l,t_k}^N)|^2$$
$$+ d\sum_{k=0}^{N-1}(\mathbb{E}[Y_{t_k}^{N,I,I} - Y_{t_k}^N]^2 - \mathbb{E}[\mathbb{E}_k(Y_{t_{k+1}}^{N,I,I} - Y_{t_{k+1}}^N)]^2).$$

Owing to (21) and (24), we obtain

$$\mathbb{E}|Y_{t_k}^{N,I,I} - Y_{t_k}^N|^2 - \mathbb{E}[\mathbb{E}_k(Y_{t_{k+1}}^{N,I,I} - Y_{t_{k+1}}^N)]^2$$
$$\leq Ch^{2I-1}\mathcal{A}^N(S_0)$$
$$+ C\mathbb{E}|\mathcal{R}_{p_{0,k}}(Y_{t_k}^N)|^2 + [(1+h)(1+\gamma h) - 1]\mathbb{E}|\mathbb{E}_k[Y_{t_{k+1}}^{N,I,I} - Y_{t_{k+1}}^N]|^2$$
$$+ Ch\left(h + \frac{1}{\gamma}\right)[\mathbb{E}|Y_{t_k}^{N,\infty,I} - Y_{t_k}^N|^2 + \mathbb{E}|Z_{t_k}^{N,I,I} - Z_{t_k}^N|^2].$$

Taking $\gamma = 4Cd$ and $h$ small enough such that $dC(h + \frac{1}{\gamma}) \leq \frac{1}{2}$, we have proved

$$\zeta^N \leq Ch^{2I-2}\mathcal{A}^N(S_0) + Ch\sum_{k=0}^{N-1}\sum_{l=1}^{q}\mathbb{E}|\mathcal{R}_{p_{l,k}}(Z_{l,t_k}^N)|^2 + C\sum_{k=0}^{N-1}\mathbb{E}|\mathcal{R}_{p_{0,k}}(Y_{t_k}^N)|^2$$
$$+ C\max_{0\leq k\leq N-1}\eta_k^{N,I} + \tfrac{1}{2}h\sum_{k=0}^{N-1}\mathbb{E}|Y_{t_k}^{N,\infty,I} - Y_{t_k}^N|^2 + \tfrac{1}{2}\zeta^N.$$

But taking into account (22) and (25) to estimate $\mathbb{E}|Y_{t_k}^{N,\infty,I} - Y_{t_k}^N|^2$, we clearly obtain (26). This easily completes the proof of Theorem 2. □

**5. Approximation results: step 3.** This step is very analogous to step 2, except that in the sequence of iterative least squares problems (9), the expectation $\mathbb{E}$ is replaced by an empirical mean built on $M$ independent simulations of $(P_{t_k}^N)_{0\leq k\leq N}, (\Delta W_k)_{0\leq k\leq N-1}$. This leads to the algorithm that is presented at Section 2.3. In this procedure, some truncation functions $\hat{\rho}_{l,k}^N$ and $\hat{\rho}_{l,k}^{N,m}$ are used and we have to specify them now.

These truncations come from a priori estimates on $Y_{t_k}^{N,i,I}, Z_{l,t_k}^{N,i,I}$ and it is useful to force their simulation-based evaluations $Y_{t_k}^{N,i,I,M,m}, Z_{l,t_k}^{N,i,I,M,m}$ to satisfy the same estimates. These a priori estimates are given by the following result (which is proved later).

PROPOSITION 2. *Under* (H1)–(H3), *for some constant $C_0$ large enough, the sequence of functions $(\rho_{l,k}^N(\cdot) = \max(1, C_0|p_{l,k}(\cdot)|) : 0 \leq l \leq q, 0 \leq k \leq N-1)$ is such that*

$$|Y_{t_k}^{N,i,I}| \leq \rho_{0,k}^N(P_{t_k}^N), \qquad \sqrt{h}|Z_{l,t_k}^{N,i,I}| \leq \rho_{l,k}^N(P_{t_k}^N) \qquad a.s.,$$

*for any $i \geq 0$, $I \geq 0$ and $0 \leq k \leq N-1$.*

With the notation of Section 2, the definition of the (random) truncation functions $\hat\rho_{l,k}^N$ (resp. $\hat\rho_{l,k}^{N,m}$) follows. Note that they are such that:

- they leave invariant $\alpha_{0,k}^{i,I} \cdot p_{0,k} = Y_{t_k}^{N,i,I}$ if $l = 0$ or $\sqrt{h}\alpha_{l,k}^{i,I} \cdot p_{l,k} = \sqrt{h}Z_{l,t_k}^{N,i,I}$ if $l \geq 1$ (resp. $\alpha_{0,k}^{i,I} \cdot p_{0,k}^m$ if $l = 0$ or $\sqrt{h}\alpha_{l,k}^{i,I} \cdot p_{l,k}^m$ if $l \geq 1$);
- they are bounded by $2\rho_{l,k}^N(P_{t_k}^N)$ [resp. $2\rho_{l,k}^N(P_{t_k}^{N,m})$];
- their first derivative is bounded by 1;
- their second derivative is uniformly bounded in $N, l, k, m$.

Now, we aim at quantifying the error between $(Y_{t_k}^{N,I,I,M}, \sqrt{h}Z_{l,t_k}^{N,I,I,M})_{l,k}$ and $(Y_{t_k}^{N,I,I}, \sqrt{h}Z_{l,t_k}^{N,I,I})_{l,k}$, in terms of the number of simulations $M$, the function bases and the time step $h$. The analysis here is more involved than in [6] since all the regression operators are correlated by the same set of simulated paths. To obtain more tractable theoretical estimates, we shall assume that each function basis $p_{l,k}$ is orthonormal. Of course, this hypothesis does not affect the numerical scheme, since the projection on a function basis is unchanged by any linear transformation of the basis. Moreover, we define the event

$$\mathbf{A}_k^M = \{\forall j \in \{k, \dots, N-1\} : \|V_j^M - \mathrm{Id}\| \leq h, \|P_{0,j}^M - \mathrm{Id}\| \leq h \quad (27) \qquad \text{and } \|P_{l,j}^M - \mathrm{Id}\| \leq 1 \text{ for } 1 \leq l \leq q\}$$

(see the notation of Section 2 for the definition of the matrices $V_j^M$ and $P_{l,j}^M$). Under the orthonormality assumption for each basis $p_{l,k}$, the matrices $(V_k^M)_{0 \leq k \leq N-1}$, $(P_{l,k}^M)_{0 \leq l \leq q, 0 \leq k \leq N-1}$ converge to the identity with probability 1 as $M \to \infty$. Thus, we have $\lim_{M \to \infty} \mathbb{P}(\mathbf{A}_k^M) = 1$. We now state our main result about the influence of the number of simulations.

THEOREM 3. *Assume* (H1)–(H3), $I \geq 3$, *that each function basis $p_{l,k}$ is orthonormal and that $\mathbb{E}|p_{l,k}|^4 < \infty$ for any $k, l$. For $h$ small enough, we have,*

*for any $0 \leq k \leq N-1$,*

$$\mathbb{E}|Y_{t_k}^{N,I,I} - Y_{t_k}^{N,I,I,M}|^2 + h\sum_{j=k}^{N-1}\mathbb{E}|Z_{t_j}^{N,I,I} - Z_{t_j}^{N,I,I,M}|^2$$

$$\leq 9\sum_{j=k}^{N-1}\mathbb{E}(|\rho_j^N(P_{t_j}^N)|^2 \mathbf{1}_{[\mathbf{A}_k^M]^c}) + Ch^{I-1}\sum_{j=k}^{N-1}[1+|S_0|^2 + \mathbb{E}|\rho_j^N(P_{t_j}^N)|^2]$$

$$+ \frac{C}{hM}\sum_{j=k}^{N-1}\bigg(\mathbb{E}\|v_j v_j^* - \mathrm{Id}\|_F^2 \mathbb{E}|\rho_j^N(P_{t_j}^N)|^2$$

$$+ \mathbb{E}(|v_j|^2|p_{0,j+1}|^2)\mathbb{E}|\rho_{0,j}^N(P_{t_j}^N)|^2$$

$$+ h^2\mathbb{E}\bigg[|v_j|^2(1+|S_{t_j}^N|^2 + |p_{0,j}|^2\mathbb{E}|\rho_{0,j}^N(P_{t_j}^N)|^2$$

$$+ \frac{1}{h}\sum_{l=1}^q |p_{l,j}|^2\mathbb{E}|\rho_{l,j}^N(P_{t_j}^N)|^2)\bigg]\bigg).$$

The term with $[\mathbf{A}_k^M]^c$ readily converges to 0 as $M \to \infty$, but we have not made estimations more explicit because the derivation of an optimal upper bound essentially depends on extra moment assumptions that may be available. For instance, if $\rho_j^N(P_{t_j}^N)$ has moments of order higher than 2, we are reduced via Hölder inequality to estimate the probability $\mathbb{P}([\mathbf{A}_k^M]^c) \leq \sum_{j=k}^{N-1}[\mathbb{P}(\|V_j^M - \mathrm{Id}\| > h) + \mathbb{P}(\|P_{0,j}^M - \mathrm{Id}\| > h) + \sum_{l=1}^q \mathbb{P}(\|P_{l,j}^M - \mathrm{Id}\| > 1)]$. We have $\mathbb{P}(\|V_k^M - \mathrm{Id}\| > h) \leq h^{-2}\mathbb{E}\|V_k^M - \mathrm{Id}\|^2 \leq h^{-2}\mathbb{E}\|V_k^M - \mathrm{Id}\|_F^2 = (Mh^2)^{-1}\mathbb{E}\|v_k v_k^* - \mathrm{Id}\|_F^2$. This simple calculus illustrates the possible computations, other terms can be handled analogously.

The previous theorem is really informative since it provides a nonasymptotic error estimation. With Theorems 1 and 2, it enables to see how to optimally choose the time step $h$, the function bases and the number of simulations to achieve a given accuracy. We do not report this analysis which seems to be hard to derive for general function bases. This will be addressed in further researches [19]. However, our next numerical experiments give an idea of this optimal choice.

We conclude our theoretical analysis by stating a central limit theorem on the coefficients $\alpha_k^{i,I,M}$ as $M$ goes to $\infty$. This is less informative than Theorem 3 since this is an asymptotic result. Thus, we remain vague about the asymptotic variance. Explicit expressions can be derived from the proof.

THEOREM 4. *Assume* (H1)–(H3), *that the driver is continuously differentiable w.r.t.* $(y,z)$ *with a bounded and uniformly Hölder continuous* derivatives and that $\mathbb{E}|p_{l,k}|^{2+\varepsilon} < \infty$ for any $k,l$ ($\varepsilon > 0$). Then, the vector $[\sqrt{M}(\alpha_k^{i,I,M} - \alpha_k^{i,I})]_{i \leq I, k \leq N-1}$ weakly converges to a centered Gaussian vector as $M$ goes to $\infty$.

PROOF OF PROPOSITION 2. In view of Proposition 1, it is tempting to apply a Markov property argument and to assert that Proposition 2 results from (19) written with conditional expectations $\mathbb{E}_k$. But this argumentation fails because the law used for the projection is not the conditional law $\mathbb{E}_k$ but $\mathbb{E}_0$. The right argument may be the following one. Write $Y_{t_k}^{N,i,I} = \alpha_{0,k}^{i,I} \cdot p_{0,k}(P_{t_k}^N)$. On the one hand, by (19), we have $C\mathcal{A}^N(S_0) \geq \mathbb{E}|Y_{t_k}^{N,i,I}|^2 = \alpha_{0,k}^{i,I} \cdot \mathbb{E}[p_{0,k}p_{0,k}^*]\alpha_{0,k}^{i,I} \geq |\alpha_{0,k}^{i,I}|^2 \lambda_{\min}(\mathbb{E}[p_{0,k}p_{0,k}^*])$. On the other hand, $|Y_{t_k}^{N,i,I}| \leq |\alpha_{0,k}^{i,I}||p_{0,k}(P_{t_k}^N)| \leq |p_{0,k}|\sqrt{C\mathcal{A}^N(S_0)/\lambda_{\min}(\mathbb{E}[p_{0,k}p_{0,k}^*])}$. Thus, we can take $\rho_{0,k}^N(x) = \max(1, |p_{0,k}(x)|\sqrt{C\mathcal{A}^N(S_0)/\lambda_{\min}(\mathbb{E}[p_{0,k}p_{0,k}^*])})$. Analogously, for $\sqrt{h}|Z_{l,t_k}^{N,i,I}|$, we have $\rho_{l,k}^N(x) = \max(1, |p_{l,k}(x)|\sqrt{C\mathcal{A}^N(S_0)/\lambda_{\min}(\mathbb{E}[p_{l,k}p_{l,k}^*])})$. Note that if $p_{l,k}$ is an orthonormal function basis, we have $\lambda_{\min}(\mathbb{E}[p_{l,k}p_{l,k}^*]) = 1$ and previous upper bounds have simpler expressions. □

PROOF OF THEOREM 3. In the sequel, set

$$\mathcal{A}_k^{N,M} = \frac{1}{M}\sum_{m=1}^M |\rho_{0,k}^N(P_{t_k}^{N,m})|^2, \qquad \mathcal{B}_k^{N,M} = \frac{1}{M}\sum_{m=1}^M |f_k^m(0,\ldots,0)|^2.$$

Obviously, we have $\mathbb{E}(\mathcal{A}_k^{N,M}) = \mathbb{E}|\rho_{0,k}^N(P_{t_k}^N)|^2$ and $\mathbb{E}(\mathcal{B}_k^{N,M}) \leq C(1 + |S_0|^2)$. Now, we remind the standard contraction property in the case of least squares problems in $\mathbb{R}^M$, analogously to the case $\mathbf{L}_2(\Omega, \mathbb{P})$. Consider a sequence of real numbers $(x^m)_{1 \leq m \leq M}$ and a sequence $(v^m)_{1 \leq m \leq M}$ of vectors in $\mathbb{R}^n$, associated to the matrix $V^M = \frac{1}{M}\sum_{m=1}^M v^m[v^m]^*$ which is supposed to be invertible $[\lambda_{\min}(V^M) > 0]$. Then, the (unique) $\mathbb{R}^n$-valued vector $\theta_x = \arg\inf_\theta |x - \theta \cdot v|_M^2$ is given by

(28) $$\theta_x = \frac{[V^M]^{-1}}{M}\sum_{m=1}^M v^m x^m.$$

The application $x \mapsto \theta_x$ is linear and, moreover, we have the inequality

(29) $$\lambda_{\min}(V^M)|\theta_x|^2 \leq |\theta_x \cdot v|_M^2 \leq |x|_M^2.$$

For the further computations, it is more convenient to deal with

$$(\theta_k^{i,I,M})^* = (\alpha_{0,k}^{i,I,M*}, \sqrt{h}\alpha_{1,k}^{i,I,M*}, \ldots, \sqrt{h}\alpha_{q,k}^{i,I,M*})$$

instead of $\alpha_k^{i,I,M}$. Then, the Picard iterations given in (4) can be rewritten

$$(30) \quad \theta_k^{i+1,I,M} = \arg\inf_\theta \frac{1}{M} \sum_{m=1}^M (\hat{\rho}_{0,k+1}^{N,m}(\alpha_{0,k+1}^{I,I,M} \cdot p_{0,k+1}^m) + h f_k^m(\alpha_k^{i,I,M}) - \theta \cdot v_k^m)^2.$$

Introducing the event $\mathbf{A}_k^M$, taking into account the Lipschitz property of the functions $\hat{\rho}_{l,k}^N$ and using the orthonormality of $p_{l,k}$, we get

$$\mathbb{E}|Y_{t_k}^{N,I,I} - Y_{t_k}^{N,I,I,M}|^2 + h \sum_{j=k}^{N-1} \mathbb{E}|Z_{t_j}^{N,I,I} - Z_{t_j}^{N,I,I,M}|^2$$

$$(31) \quad \leq 9 \sum_{j=k}^{N-1} \mathbb{E}(|\rho_j^N(P_{t_j}^N)|^2 \mathbf{1}_{[\mathbf{A}_k^M]^c})$$

$$+ \mathbb{E}(\mathbf{1}_{\mathbf{A}_k^M}|\alpha_{0,k}^{I,I,M} - \alpha_{0,k}^{I,I}|^2) + h \sum_{j=k}^{N-1} \sum_{l=1}^q \mathbb{E}(\mathbf{1}_{\mathbf{A}_k^M}|\alpha_{l,j}^{I,I,M} - \alpha_{l,j}^{I,I}|^2).$$

To obtain Theorem 3, we estimate $|\theta_k^{I,I,M} - \theta_k^{I,I}|^2$ on the event $\mathbf{A}_k^M$. This is achieved in several steps.

*Step* 1: *contraction properties relative to the sequence* $(\theta_k^{i,I,M})_{i\geq 0}$. They are summed up in the following lemma:

LEMMA 1. *For $h$ small enough, on $\mathbf{A}_k^M$ the following properties hold:*
(a) $|\theta_k^{i+1,I,M} - \theta_k^{i,I,M}|^2 \leq Ch|\theta_k^{i,I,M} - \theta_k^{i-1,I,M}|^2.$
(b) *There is a unique vector $\theta_k^{\infty,I,M}$ such that*

$$\theta_k^{\infty,I,M} = \arg\inf_\theta \frac{1}{M} \sum_{m=1}^M (\hat{\rho}_{0,k+1}^{N,m}(\alpha_{0,k+1}^{I,I,M} \cdot p_{0,k+1}^m) + h f_k^m(\alpha_k^{\infty,I,M}) - \theta \cdot v_k^m)^2.$$

(c) *We have $|\theta_k^{\infty,I,M} - \theta_k^{I,I,M}|^2 \leq [Ch]^I |\theta_k^{\infty,I,M}|^2.$*

PROOF. We prove (a). Since $1 - h \leq \lambda_{\min}(V_k^M)$ and $\lambda_{\max}(P_{l,k}^M) \leq 2$ ($0 \leq l \leq q$) on $\mathbf{A}_k^M$, in view of (29), we obtain that $(1-h)|\theta_k^{i+1,I,M} - \theta_k^{i,I,M}|^2$ is bounded by

$$\frac{h^2}{M} \sum_{m=1}^M (f_k^m(\alpha_k^{i,I,M}) - f_k^m(\alpha_k^{i-1,I,M}))^2$$

$$\leq Ch^2 \sum_{l=0}^q |\alpha_{l,k}^{i,I,M} - \alpha_{l,k}^{i-1,I,M}|^2 \lambda_{\max}(P_{l,k}^M)$$

$$\leq Ch|\theta_k^{i,I,M} - \theta_k^{i-1,I,M}|^2.$$

Now, statements (a) and (b) are clear. For (c), apply (a), reminding that $\theta_k^{0,I,M} = 0$. □

*Step* 2: *bounds for* $|\theta_k^{i,I,M}|$ *on the event* $\mathbf{A}_k^M$. Namely, we aim at showing that

$$|\theta_k^{i,I,M}|^2 \leq C(\mathcal{A}_{k+1}^{N,M} + h\mathcal{B}_k^{N,M}) \qquad \text{on } \mathbf{A}_k^M. \tag{32}$$

We first consider $i = \infty$. As in the proof of Lemma 1, we get

$$(1-h)|\theta_k^{\infty,I,M}|^2$$

$$\leq \frac{1}{M}\sum_{m=1}^{M}[\hat{\rho}_{0,k+1}^{N,m}(\alpha_{0,k+1}^{I,I,M} \cdot p_{0,k+1}^m) + hf_k^m(\alpha_k^{\infty,I,M})]^2$$

$$\leq (1+\gamma h)\mathcal{A}_{k+1}^{N,M} + Ch\left(h + \frac{1}{\gamma}\right)\left(\mathcal{B}_k^{N,M} + \sum_{l=0}^{q}|\alpha_{l,k}^{\infty,I,M}|^2\lambda_{\max}(P_{l,k}^M)\right).$$

Take $\gamma = 8C$ and $h$ small enough to ensure $2C(h + \frac{1}{\gamma})(1+h) \leq \frac{1}{2}(1-h)$. It readily follows $|\theta_k^{\infty,I,M}|^2 \leq C(\mathcal{A}_{k+1}^{N,M} + h\mathcal{B}_k^{N,M})$, proving that (32) holds for $i = \infty$. Lemma 1(c) leads to expected bounds for other values of $i$.

*Step* 3: *we remind bounds for* $\theta^{i,I}$. Using Proposition 2 and in view of (10)–(14), we have, for $i \geq 1$,

$$|\theta_{l,k}^{i,I}|^2 \leq \mathbb{E}|\rho_{l,k}^N(P_{t_k}^N)|^2, \qquad 0 \leq l \leq q;$$

$$|\theta_k^{\infty,I} - \theta_k^{i,I}|^2 \leq (C_f h)^{2i}\mathbb{E}|\rho_{0,k}^N(P_{t_k}^N)|^2. \tag{33}$$

Remember also the following expression of $\theta_k^{\infty,I}$, derived from (10)–(13) and the orthonormality of each basis $p_{l,k}$:

$$\theta_k^{\infty,I} = \mathbb{E}(v_k[\alpha_{0,k+1}^{I,I} \cdot p_{0,k+1} + hf_k(\alpha_k^{\infty,I})]). \tag{34}$$

*Step* 4: *decomposition of the quantity* $\mathbb{E}(\mathbf{1}_{\mathbf{A}_k^M}|\theta_k^{I,I,M} - \theta_k^{I,I}|^2)$. Due to Lemma 1, on $\mathbf{A}_k^M$ we get $|\theta_k^{\infty,I,M} - \theta_k^{I,I,M}|^2 \leq Ch^I|\theta_k^{\infty,I,M}|^2 \leq Ch^I|\theta_k^{\infty,I}|^2 + Ch^I|\theta_k^{\infty,I,M} - \theta_k^{\infty,I}|^2$. Thus, using (33), it readily follows that $\mathbb{E}(\mathbf{1}_{\mathbf{A}_k^M}|\theta_k^{I,I,M} - \theta_k^{I,I}|^2)$ is bounded by

$$(1+h)\mathbb{E}(\mathbf{1}_{\mathbf{A}_k^M}|\theta_k^{\infty,I,M} - \theta_k^{\infty,I}|^2)$$

$$+ 2\left(1 + \frac{1}{h}\right)\{\mathbb{E}(\mathbf{1}_{\mathbf{A}_k^M}|\theta_k^{I,I,M} - \theta_k^{\infty,I,M}|^2) + |\theta_k^{I,I} - \theta_k^{\infty,I}|^2\} \tag{35}$$

$$\leq (1+Ch)\mathbb{E}(\mathbf{1}_{\mathbf{A}_k^M}|\theta_k^{\infty,I,M} - \theta_k^{\infty,I}|^2) + Ch^{I-1}\mathbb{E}|\rho_k^N(P_{t_k}^N)|^2,$$

taking account that $I \geq 3$. On $\mathbf{A}_k^M$, $V_k^M$ is invertible and we can set
$B_1 = (\mathrm{Id} - (V_k^M)^{-1})\theta_k^{\infty,I}$,

$$B_2 = (V_k^M)^{-1}\left[\mathbb{E}(v_k\hat{\rho}_{0,k+1}^N(\alpha_{0,k+1}^{I,I} \cdot p_{0,k+1})) - \frac{1}{M}\sum_{m=1}^M v_k^m \hat{\rho}_{0,k+1}^{N,m}(\alpha_{0,k+1}^{I,I} \cdot p_{0,k+1}^m)\right],$$

$$B_3 = (V_k^M)^{-1}h\left[\mathbb{E}(v_k f_k(\alpha_k^{\infty,I})) - \frac{1}{M}\sum_{m=1}^M v_k^m f_k^m(\alpha_k^{\infty,I})\right],$$

$$B_4 = \frac{(V_k^M)^{-1}}{M}\sum_{m=1}^M v_k^m[\hat{\rho}_{0,k+1}^{N,m}(\alpha_{0,k+1}^{I,I} \cdot p_{0,k+1}^m) - \hat{\rho}_{0,k+1}^{N,m}(\alpha_{0,k+1}^{I,I,M} \cdot p_{0,k+1}^m)$$

$$+ h(f_k^m(\alpha_k^{\infty,I}) - f_k^m(\alpha_k^{\infty,I,M}))].$$

Thus, by (28)–(34), we can write $\theta_k^{\infty,I} - \theta_k^{\infty,I,M} = B_1 + B_2 + B_3 + B_4$, which gives on $\mathbf{A}_k^M$

$$(36) \quad |\theta_k^{\infty,I} - \theta_k^{\infty,I,M}|^2 \leq 3\left(1 + \frac{1}{h}\right)(|B_1|^2 + |B_2|^2 + |B_3|^2) + (1+h)|B_4|^2.$$

*Step* 5: *individual estimation of* $B_1$, $B_2$, $B_3$, $B_4$ *on* $\mathbf{A}_k^M$. Remember the classic result [16]: if $\|\mathrm{Id} - F\| < 1$, $F^{-1} - \mathrm{Id} = \sum_{k=1}^\infty [\mathrm{Id} - F]^k$ and $\|\mathrm{Id} - F^{-1}\| \leq \frac{\|F-\mathrm{Id}\|}{1-\|F-\mathrm{Id}\|}$. Consequently, for $F = V_k^M$, we get $\mathbb{E}(\mathbf{1}_{\mathbf{A}_k^M}\|\mathrm{Id}-(V_k^M)^{-1}\|^2) \leq (1-h)^{-2}\mathbb{E}\|\mathrm{Id} - V_k^M\|^2 \leq (1-h)^{-2}\mathbb{E}\|V_k^M - \mathrm{Id}\|_F^2 = (M(1-h)^2)^{-1}\mathbb{E}\|v_kv_k^* - \mathrm{Id}\|_F^2$. Thus, we have

$$\mathbb{E}(|B_1|^2 \mathbf{1}_{\mathbf{A}_k^M}) \leq \frac{C}{M}\mathbb{E}\|v_kv_k^* - \mathrm{Id}\|_F^2 \mathbb{E}|\rho_k^N(P_{t_k}^N)|^2.$$

Since on $\mathbf{A}_k^M$ one has $\|(V_k^M)^{-1}\| \leq 2$, it readily follows

$$\mathbb{E}(|B_2|^2 \mathbf{1}_{\mathbf{A}_k^M}) \leq \frac{C}{M}\mathbb{E}(|v_k|^2|p_{0,k+1}|^2)\mathbb{E}|\rho_{0,k}^N(P_{t_k}^N)|^2,$$

$$\mathbb{E}(|B_3|^2 \mathbf{1}_{\mathbf{A}_k^M}) \leq \frac{Ch^2}{M}\mathbb{E}\left[|v_k|^2\left(1 + |S_{t_k}^N|^2 + |p_{0,k}|^2\mathbb{E}|\rho_{0,k}^N(P_{t_k}^N)|^2\right.\right.$$

$$\left.\left.+ \frac{1}{h}\sum_{l=1}^q |p_{l,k}|^2 \mathbb{E}|\rho_{l,k}^N(P_{t_k}^N)|^2\right)\right].$$

As in the proof of Lemma 1 and using $\|P_{0,k+1}^M\| \leq 1 + h$ on $\mathbf{A}_k^M$, we easily obtain

$$(1-h)|B_4|^2 \leq (1+h)(1+\gamma h)|\alpha_{0,k+1}^{I,I} - \alpha_{0,k+1}^{I,I,M}|^2$$

$$+ Ch\left(h + \frac{1}{\gamma}\right)\sum_{l=0}^q |\alpha_{l,k}^{\infty,I} - \alpha_{l,k}^{\infty,I,M}|^2.$$

*Step* 6: *final estimations.* Put $\epsilon_k = \mathbb{E}\|v_k v_k^* - \mathrm{Id}\|_F^2 \mathbb{E}|\rho_k^N(P_{t_k}^N)|^2 + \mathbb{E}(|v_k|^2|p_{0,k+1}|^2)\mathbb{E}|\rho_{0,k}^N(P_{t_k}^N)|^2 + h^2 \mathbb{E}[|v_k|^2(1+|S_{t_k}^N|^2+|p_{0,k}|^2\mathbb{E}|\rho_{0,k}^N(P_{t_k}^N)|^2 + \frac{1}{h}\sum_{l=1}^{q}|p_{l,k}|^2\mathbb{E}|\rho_{l,k}^N(P_{t_k}^N)|^2)]$. Plug the above estimates on $B_1, B_2, B_3, B_4$ into (36), choose $\gamma = 3C$ and $h$ close to 0 to ensure $Ch + \frac{C}{\gamma} \leq \frac{1}{2}$; after simplifications, we get

$$\mathbb{E}(\mathbf{1}_{\mathbf{A}_k^M}|\theta_k^{\infty,I,M} - \theta_k^{\infty,I}|^2) \leq C\frac{\epsilon_k}{hM} + (1+Ch)\mathbb{E}(\mathbf{1}_{\mathbf{A}_k^M}|\alpha_{0,k+1}^{I,I} - \alpha_{0,k+1}^{I,I,M}|^2).$$

But in view of Lemma 1(c) and estimates (32)–(33), we have

$$\mathbb{E}(\mathbf{1}_{\mathbf{A}_k^M}|\alpha_{0,k+1}^{I,I} - \alpha_{0,k+1}^{I,I,M}|^2)$$
$$\leq (1+h)\mathbb{E}(\mathbf{1}_{\mathbf{A}_k^M}|\alpha_{0,k+1}^{\infty,I} - \alpha_{0,k+1}^{\infty,I,M}|^2)$$
$$+ Ch^{I-1}(1+|S_0|^2 + \mathbb{E}|\rho_{0,k+1}^N(P_{t_{k+1}}^N)|^2 + \mathbb{E}|\rho_{0,k+2}^N(P_{t_{k+2}}^N)|^2).$$

Finally, we have proved

$$\mathbb{E}(\mathbf{1}_{\mathbf{A}_k^M}|\theta_k^{\infty,I,M} - \theta_k^{\infty,I}|^2)$$
$$\leq C\frac{\epsilon_k}{hM} + Ch^{I-1}(1+|S_0|^2 + \mathbb{E}|\rho_{0,k+1}^N(P_{t_{k+1}}^N)|^2 + \mathbb{E}|\rho_{0,k+2}^N(P_{t_{k+2}}^N)|^2)$$
$$+ (1+Ch)\mathbb{E}(\mathbf{1}_{\mathbf{A}_k^M}|\alpha_{0,k+1}^{\infty,I,M} - \alpha_{0,k+1}^{\infty,I}|^2).$$

Using a contraction argument as in (35), the index $\infty$ can be replaced by $I$, without changing the inequality (with a possibly different constant $C$). This can be written

$$\mathbb{E}(\mathbf{1}_{\mathbf{A}_k^M}|\alpha_{0,k}^{I,I,M} - \alpha_{0,k}^{I,I}|^2) + h\sum_{l=1}^{q}\mathbb{E}(\mathbf{1}_{\mathbf{A}_k^M}|\alpha_{l,k}^{I,I,M} - \alpha_{l,k}^{I,I}|^2)$$
$$\leq C\frac{\epsilon_k}{hM} + Ch^{I-1}(1+|S_0|^2 + \mathbb{E}|\rho_{0,k+1}^N(P_{t_{k+1}}^N)|^2 + \mathbb{E}|\rho_{0,k+2}^N(P_{t_{k+2}}^N)|^2)$$
$$+ (1+Ch)\mathbb{E}(\mathbf{1}_{\mathbf{A}_k^M}|\alpha_{0,k+1}^{I,I,M} - \alpha_{0,k+1}^{I,I}|^2).$$

Using Gronwall's lemma, the proof is complete. □

REMARK 1. The attentive reader may have noted that powers of $h$ are smaller here than in Theorem 2, which leads to take $I \geq 3$ instead of $I \geq 2$ before. Indeed, we cannot take advantage of conditional expectations on the simulations as we did in (12), for instance.

Note that in the proof above, we only use the Lipschitz property of the truncation functions $\hat{\rho}_{l,k}^N$ and $\hat{\rho}_{l,k}^{N,m}$.

PROOF OF THEOREM 4. The arguments are standard and there are essentially notational difficulties. The first partial derivatives of $f$ w.r.t. $y$

and $z_l$ are, respectively, denoted $\partial_0 f$ and $\partial_l f$. The parameter $\beta \in ]0,1]$ stands for their Hölder continuity index. Suppose w.l.o.g. that $\varepsilon < \beta$ and that each function basis $p_{l,k}$ is orthonormal. For $k < N-1$, define the quantities

$$A_{l,k}^M(\alpha) = \frac{1}{M} \sum_{m=1}^{M} v_k^m \, \partial_l f(t_k, S_{t_k}^{N,m}, \alpha_0 \cdot p_{0,k}^m, \ldots, \alpha_q \cdot p_{q,k}^m)[p_{l,k}^m]^*,$$

$$B_k^M = \frac{1}{M} \sum_{m=1}^{M} v_k^m [p_{0,k+1}^m]^*, \qquad D_k^M = \sqrt{M}(\mathrm{Id} - V_k^M),$$

$$C_k^M(\alpha)$$
$$= \sum_{m=1}^{M} \frac{\{v_k^m [\alpha_{0,k+1}^{I,I} \cdot p_{0,k+1}^m + h f_k^m(\alpha)] - \mathbb{E}(v_k[\alpha_{0,k+1}^{I,I} \cdot p_{0,k+1} + h f_k(\alpha)])\}}{\sqrt{M}}.$$

For $k = N-1$, we set $B_k^M = 0$ and in $C_k^M(\alpha)$, the terms $\alpha_{0,k+1}^{I,I} \cdot p_{0,k+1}^m$ and $\alpha_{0,k+1}^{I,I} \cdot p_{0,k+1}$ have to be replaced, respectively, by $\Phi^N(P_{t_N}^{N,m})$ and $\Phi^N(P_{t_N}^N)$. The definitions of $A_{l,k}^M(\alpha)$ and $D_k^M$ are still valid. For convenience, we write $X^M \xrightarrow{w}$ if the (possibly vector or matrix valued) sequence $(X^M)_M$ weakly converges to a centered Gaussian variable, as $M$ goes to infinity. For the convergence in probability to a constant, we denote $X^M \xrightarrow{\mathbb{P}}$. Since simulations are independent, observe that the following convergences hold:

(37)
$$(A_{l,k}^M(\alpha_k^{i,I}), B_k^M, V_k^M)_{i \leq I-1, l \leq q, k \leq N-1} \xrightarrow{\mathbb{P}},$$
$$\mathcal{G}^M = (C_k^M(\alpha_k^{i,I}), D_k^M)_{i \leq I-1, l \leq q, k \leq N-1} \xrightarrow{w}.$$

Note that $\lim_{M \to \infty} V_k^M \stackrel{\text{a.s.}}{=} \mathrm{Id}$ is invertible. Linearizing the functions $f$ and $\hat{\rho}_{0,k+1}^{N,m}$ in the expressions of $\theta_k^{i,I} = \mathbb{E}(v_k[\alpha_{0,k+1}^{I,I} \cdot p_{0,k+1} + h f_k(\alpha_{0,k}^{i-1,I}, \ldots, \alpha_{q,k}^{i-1,I})])$ and $\theta_k^{i,I,M}$ given by (28) leads to

(38)
$$\left| V_k^M \sqrt{M}(\theta_k^{i,I,M} - \theta_k^{i,I}) - D_k^M \theta_k^{i,I} - C_k^M(\alpha_k^{i-1,I}) \right.$$
$$- B_k^M \sqrt{M}(\alpha_{0,k+1}^{I,I,M} - \alpha_{0,k+1}^{I,I})$$
$$\left. - h \sum_{l=0}^{q} A_{l,k}^M(\alpha_k^{i-1,I}) \sqrt{M}(\alpha_{l,k}^{i-1,I,M} - \alpha_{l,k}^{i-1,I}) \right|$$
$$\leq \mathbf{1}_{k<N-1} \frac{C}{\sqrt{M}} |\alpha_{0,k+1}^{I,I,M} - \alpha_{0,k+1}^{I,I}|^2 \sum_{m=1}^{M} |v_k^m||p_{0,k+1}^m|^2$$

$$+ \frac{C}{\sqrt{M}}|\alpha_k^{i-1,I,M} - \alpha_k^{i-1,I}|^{1+\beta} \sum_{m=1}^M |v_k^m||p_k^m|^{1+\beta}.$$

To get Theorem 4, we prove by induction on $k$ that $([\sqrt{M}(\theta_j^{i,I,M} - \theta_j^{i,I})]_{j \geq k, i \leq I}, \mathcal{G}^M) \xrightarrow{w}$. Remember that $\theta_j^{0,I,M} = \theta_j^{0,I} = 0$ for any $j$. Consider first $k = N-1$, for which $B_k^M = 0$, and $i = 1$. In view of (37)–(38), clearly $([\sqrt{M}(\theta_{N-1}^{i,I,M} - \theta_{N-1}^{i,I})]_{i \leq 1}, \mathcal{G}^M) \xrightarrow{w}$. For $i = 2$, we may invoke the same argument using (37)–(38) and obtain $([\sqrt{M}(\theta_{N-1}^{i,I,M} - \theta_{N-1}^{i,I})]_{i \leq 2}, \mathcal{G}^M) \xrightarrow{w}$ provided that the upper bound in (38) converge to 0 in probability. To prove this, put $\mathcal{M}^M = M^{-1-\beta/2} \times \sum_{m=1}^M |v_{N-1}^m||p_{N-1}^m|^{1+\beta}$ and write $\frac{1}{\sqrt{M}}|\alpha_{N-1}^{1,I,M} - \alpha_{N-1}^{1,I}|^{1+\beta} \sum_{m=1}^M |v_{N-1}^m| \times |p_{N-1}^m|^{1+\beta} = |\sqrt{M}(\alpha_{N-1}^{1,I,M} - \alpha_{N-1}^{1,I})|^{1+\beta} \mathcal{M}^M$. Since $[\sqrt{M}(\alpha_{N-1}^{1,I,M} - \alpha_{N-1}^{1,I})]_M$ is tight, our assertion holds if $\mathcal{M}^M$ converges to 0 as $M \to \infty$. Note that $|v_{N-1}||p_{N-1}|^{1+\beta} \in \mathbf{L}_{(2+\varepsilon)/(2+\beta)}(\mathbb{P})$. Thus, the strong law of large numbers, in the case of i.i.d. random variables with infinite mean, leads to $\sum_{m=1}^M |v_{N-1}^m| \times |p_{N-1}^m|^{1+\beta} = O(M^{(2+\beta)/(2+\varepsilon)+r})$ a.s. for any $r > 0$. Consequently, from the choice of $r$ small enough, it follows $\mathcal{M}^M \to 0$ a.s.

Iterating this argumentation readily leads to $([\sqrt{M}(\theta_{N-1}^{i,I,M} - \theta_{N-1}^{i,I})]_{i \leq \mathbf{I}}, \mathcal{G}^M) \xrightarrow{w}$. For the induction for $k < N-1$, we apply the techniques above. There is an additional contribution due to $B_k^M$, which can be handled as before. $\square$

## 6. Numerical experiments.

6.1. *Lipschitz property of the solution under* (H4). To use the algorithm, we need to specify the basis functions that we choose at each time $t_k$ and for this, the knowledge of the regularity of the functions $y_k^N(\cdot)$ and $z_{l,k}^N(\cdot)$ from Proposition 1 is useful (in view of Theorem 2). In all the cases described in Section 2.4 and below, assumption (H4) is fulfilled. Under this extra assumption, we now establish that $y_k^N(\cdot)$ and $z_{l,k}^N(\cdot)$ are Lipschitz continuous.

PROPOSITION 3. *Assume* (H1)–(H4). *For $h$ small enough, we have*

$$|y_{k_0}^N(x) - y_{k_0}^N(x')| + \sqrt{h}|z_{k_0}^N(x) - z_{k_0}^N(x')| \leq C|x - x'| \tag{39}$$

*uniformly in $k_0 \leq N - 1$.*

PROOF. As for (17), we can obtain

$$\mathbb{E}|Y_{t_k}^{N,k_0,x} - Y_{t_k}^{N,k_0,x'}|^2$$
$$\leq \frac{(1+\gamma h)}{1 - Ch(h+1/\gamma)}\mathbb{E}|\mathbb{E}_k(Y_{t_{k+1}}^{N,k_0,x} - Y_{t_{k+1}}^{N,k_0,x'})|^2$$
$$+ \frac{Ch(h+1/\gamma)}{1 - Ch(h+1/\gamma)}\mathbb{E}|S_{t_k}^{N,k_0,x} - S_{t_k}^{N,k_0,x'}|^2$$

$$+ \frac{C(h+1/\gamma)}{1-Ch(h+1/\gamma)}(\mathbb{E}|Y_{t_{k+1}}^{N,k_0,x} - Y_{t_{k+1}}^{N,k_0,x'}|^2$$
$$- \mathbb{E}|\mathbb{E}_k(Y_{t_{k+1}}^{N,k_0,x} - Y_{t_{k+1}}^{N,k_0,x'})|^2).$$

Choosing $\gamma = C$ and $h$ small enough, we get (for another constant $C$)

$$\mathbb{E}|Y_{t_k}^{N,k_0,x} - Y_{t_k}^{N,k_0,x'}|^2$$
$$\leq (1+Ch)\mathbb{E}|Y_{t_{k+1}}^{N,k_0,x} - Y_{t_{k+1}}^{N,k_0,x'}|^2 + Ch\mathbb{E}|S_{t_k}^{N,k_0,x} - S_{t_k}^{N,k_0,x'}|^2.$$

The last term above is bounded by $C|x-x'|^2$ under assumption (H1). Thus, using Gronwall's lemma and assumption (H4), we get the result for $y_{k_0}^N(\cdot)$. The result for $\sqrt{h}z_{k_0}^N(\cdot)$ follows by considering (5). □

6.2. *Choice of function bases.* Now, we specify several choices of function bases. We denote $d'$ ($\geq d$) the dimension of the state space of $(P_{t_k}^N)_k$.

*Hypercubes* (HC in the following). Here, to simplify, $p_{l,k}$ does not depend on $l$ or $k$. Choose a domain $D \subset \mathbb{R}^{d'}$ centered on $P_0^N$, that is, $D = \prod_{i=1}^{d'} ]P_{0,i}^N - R, P_{0,i}^N + R]$, and partition it into small hypercubes of edge $\delta$. Thus, $D = \bigcup_{i_1,\ldots,i_{d'}} D_{i_1,\ldots,i_{d'}}$ where $D_{i_1,\ldots,i_{d'}} = ]P_{0,1}^N - R + i_1\delta, P_{0,1}^N - R + (i_1+1)\delta] \times \cdots \times ]P_{0,d'}^N - R + i_{d'}\delta, P_{0,d'}^N - R + (i_{d'}+1)\delta]$. Then we define $p_{l,k}$ as the indicator functions associated to this set of hypercubes: $p_{l,k}(\cdot) = (\mathbf{1}_{D_{i_1,\ldots,i_{d'}}}(\cdot))_{i_1,\ldots,i_{d'}}$. With this particular choice of function bases, we can explicit the projection error of Theorem 2:

$$\mathbb{E}(\mathcal{R}_{p_{0,k}}(Y_{t_k}^N)^2)$$
$$\leq \mathbb{E}(|Y_{t_k}^N|^2 \mathbf{1}_{D^c}(P_{t_k}^N)) + \sum_{i_1,\ldots,i_{d'}} \mathbb{E}(\mathbf{1}_{D_{i_1,\ldots,i_{d'}}}(P_{t_k}^N)|y_k^N(P_{t_k}^N) - y_k^N(x_{i_1,\ldots,i_{d'}})|^2)$$
$$\leq C\delta^2 + \mathbb{E}(|Y_{t_k}^N|^2 \mathbf{1}_{D^c}(P_{t_k}^N)),$$

where $x_{i_1,\ldots,i_{d'}}$ is an arbitrary point of $D_{i_1,\ldots,i_{d'}}$ and where we have used the Lipschitz property of $y_k^N$ on $D$. To evaluate $\mathbb{E}(|Y_{t_k}^N|^2 \mathbf{1}_{D^c}(P_{t_k}^N))$, note that, by adapting the proof of Proposition 3, we have $|Y_{t_k}^N|^2 \leq C(1 + |S_{t_k}^N|^2 + \mathbb{E}_k|P_{t_N}^N|^2)$. Thus, if $\sup_{k,N} \mathbb{E}|P_{t_k}^N|^\alpha < \infty$ for $\alpha > 2$, we have $\mathbb{E}(|Y_k^N|^2 \mathbf{1}_{D^c}(P_{t_k}^N)) \leq \frac{C_\alpha}{R^{\alpha-2}}$, with an explicit constant $C_\alpha$. The choice $R \approx h^{-2/(\alpha-2)}$ and $\delta = h$ leads to

$$\mathbb{E}|\mathcal{R}_{p_{0,k}}(Y_{t_k}^N)|^2 \leq Ch^2.$$

The same estimates hold for $\mathbb{E}|\mathcal{R}_{p_{l,k}}(\sqrt{h}Z_{l,t_k}^N)|^2$. Thus, we obtain the same accuracy as in Theorem 1.

*Voronoi partition* (VP). Here, we consider again a basis of indicator functions and the same basis for all $0 \leq l \leq q$. This time, the sets of the indicator functions are an open Voronoi partition ([17]) whose centers are independent simulations of $P^N$. More precisely, if we want a basis of 20 indicator functions, we simulate 20 extra paths of $P^N$, denoted $(P^{N,M+i})_{1\leq i\leq 20}$, independently of $(P^{N,m})_{1\leq m\leq M}$. Then we take at time $t_k$ $(P_{t_k}^{N,M+i})_{1\leq i\leq 20}$ to define our Voronoi partition $(C_{k,i})_{1\leq i\leq 20}$, where $C_{k,i} = \{x : |x - P_{t_k}^{N,M+i}| < \inf_{j\neq i} |x - P_{t_k}^{N,M+j}|\}$. Then $p_{l,k}(\cdot) = (\mathbf{1}_{C_{k,i}}(\cdot))_i$. We can notice that, unlike the hypercubes basis, the function basis changes with $k$. We can also estimate the projection error of Theorem 2, using results on random quantization and refer to [17] for explicit calculations.

In addition, we can consider on each Voronoi cells local polynomials of low degree. For example, we can take a local polynomial basis consisting of $1, x_1, \ldots, x_{d'}$ for $p_{0,k}$ and 1 for $p_{l,k}$ ($l \geq 1$) on each $C_{k,i}$. Thus, $p_{0,k}(x) = (\mathbf{1}_{C_{k,i}}(x), x_1 \mathbf{1}_{C_{k,i}}(x), \ldots, x_{d'} \mathbf{1}_{C_{k,i}}(x))_i$ and $p_{l,k}(x) = (\mathbf{1}_{C_{k,i}}(x))_i, 1 \leq l \leq q$. We denote this particular choice $\mathrm{VP}(1,0)$, where 1 (resp. 0) stands for the maximal degree of local polynomial basis for $p_{0,k}$ (resp. $p_{l,k}$, $1 \leq l \leq q$).

*Global polynomials* (GP). Here we define $p_{0,k}$ as the polynomial (of $d'$ variables) basis of degree less than $d_y$ and $p_{l,k}$ as the polynomial basis of degree less than $d_z$.

6.3. *Numerical results.* After the description of possible basis functions, we test the algorithm on several examples. For each example and each choice of function basis, we launch the algorithm for different values of $M$, the number of Monte Carlo simulations. More precisely, for each value of $M$, we launch 50 times the algorithm and collect each time the value $Y_{t_0}^{N,I,I,M}$. The set of collected values is denoted $(Y_{t_0,i}^{N,I,I,M})_{1\leq i\leq 50}$. Then, we compute the empirical mean $\overline{Y}_{t_0}^{N,I,I,M} = \frac{1}{50}\sum_{i=1}^{50} Y_{t_0,i}^{N,I,I,M}$ and the empirical standard deviation $\sigma_{t_0}^{N,I,I,M} = \sqrt{\frac{1}{49}\sum_{i=1}^{50} |Y_{t_0,i}^{N,I,I,M} - \overline{Y}_{t_0}^{N,I,I,M}|^2}$. These two statistics provide an insight into the accuracy of the method.

6.3.1. *Call option with different interest rates 4.* We follow the notation of Section 2.4 considering a one-dimensional Black–Scholes model, with parameters

| $\mu$ | $\sigma$ | $r$ | $R$ | $T$ | $S_0$ | $K$ |
|---|---|---|---|---|---|---|
| 0.06 | 0.2 | 0.04 | 0.06 | 0.5 | 100 | 100 |

Here $K$ is the strike of the call option: $\Phi(\mathbf{S}) = (S_T - K)_+$. We know by the comparison theorem for BSDEs [12] and properties of the price and replicating strategies of a call option, that the seller of the option has always to borrow money to replicate the option in continuous time. Thus, $Y_0$ is given by the Black–Scholes formula evaluated with interest rate $R : \mathbf{Y}_0 = 7.15$. This is a good test for our algorithm because the driver $f$ is nonlinear, but we nevertheless know the true value of $Y_0$. We test the function basis HC for different values of $N$, $D$ and $\delta$. Results ($\overline{Y}_{t_0}^{N,I,I,M}$ and $\sigma_{t_0}^{N,I,I,M}$ in parenthesis) are reported in Table 1, for different values of $M$. Clearly, $\overline{Y}_{t_0}^{N,I,I,M}$ converges toward 7.15, which is exactly the Black–Scholes price $Y_0$ calculated with interest rate $R$. We observe a decrease of the empirical standard deviation like $1/\sqrt{M}$, which is coherent with Theorem 4.

6.3.2. *Calls combination with different interest rates.* We take the same driver $f$ but change the terminal condition: $\Phi(\mathbf{S}) = (S_T - K_1)^+ - 2(S_T - K_2)^+$. We take the following values for the parameters:

| $\mu$ | $\sigma$ | $r$ | $R$ | $T$ | $S_0$ | $K_1$ | $K_2$ |
|---|---|---|---|---|---|---|---|
| 0.05 | 0.2 | 0.01 | 0.06 | 0.25 | 100 | 95 | 105 |

We denote by $BS_i(r)$ the Black–Scholes price evaluated with strike $K_i$ and interest rate $r$. If we try to predict $Y_0$ by a linear combination of Black–Scholes prices, we get

| | |
|---|---|
| $BS_1(R) - 2BS_2(R)$ | 2.75 |
| $BS_1(r) - 2BS_2(r)$ | 2.76 |
| $BS_1(r) - 2BS_2(R)$ | 1.92 |
| $BS_1(R) - 2BS_2(r)$ | 3.60 |

Using comparison results, one can check that the first three rows provide a lower bound for $Y_0$, while the fourth row gives an upper bound. According to the results of HC and VP, $\overline{Y}_{t_0}^{N,I,I,M}$ seems to converge toward 2.95. This value is not predicted by a linear combination of Black–Scholes prices: in this example, the nonlinearity of $f$ has a real impact on $Y_0$. The financial interpretation is that the option seller has alternatively to borrow and to lend money to replicate the option payoff.

Comparing the different choices of basis functions, we can notice that the column $N = 5$ of VP (Table 3) shows similar results with an equal number of basis functions than the column $N = 5$ of HC (Table 2). In Table 3, the last two columns show that using a local polynomial basis may significantly increase the accuracy. We also remark by considering the rows $M = 128, 512$ of Table 2 that the standard deviation increases with $N$ and the number of basis functions, which is coherent with Theorem 3. Finally, from Table 4 the basis GP also succeeds in reaching the expected value, as we increase the number of polynomials in the basis.

TABLE 1
*Results for the call option using the basis HC*

| $M$ | $N = 5$, $D = [60, 140]$, $\delta = 5$ | $N = 10$, $D = [60, 140]$, $\delta = 1$ |
|---|---|---|
| 128 | 6.83(0.31) | 7.02(0.51) |
| 512 | 7.08(0.11) | 7.12(0.21) |
| 2048 | 7.13(0.05) | 7.14(0.07) |
| 8192 | 7.15(0.03) | 7.15(0.03) |
| 32768 | 7.15(0.01) | 7.15(0.02) |

TABLE 2
*Results for the calls combination using the basis HC*

| $M$ | $N = 5$ $D = [60, 140]$ $\delta = 5$ | $N = 20$ $D = [60, 200]$ $\delta = 1$ | $N = 50$ $D = [40, 200]$ $\delta = 0.5$ |
|---|---|---|---|
| 128 | 3.05(0.27) | 3.71(0.95) | 3.69(4.15) |
| 512 | 2.93(0.11) | 3.14(0.16) | 3.48(0.54) |
| 2048 | 2.92(0.05) | 3.00(0.03) | 3.08(0.12) |
| 8192 | 2.91(0.03) | 2.96(0.02) | 2.99(0.02) |
| 32768 | 2.90(0.01) | 2.95(0.01) | 2.96(0.01) |

TABLE 3
*Results for the calls combination using the bases VP and* $\mathrm{VP}(1,0)$

| $M$ | Basis VP 16 Voronoi regions $N = 5$ | Basis VP 64 Vor. reg. $N = 20$ | Basis VP 10 Vor. reg. $N = 20$ | Basis VP(1,0) 10 Vor. reg. $N = 20$ |
|---|---|---|---|---|
| 128 | 3.23(0.30) | 4.50(1.71) | 3.08(0.25) | 3.23(0.23) |
| 512 | 3.05(0.13) | 3.36(0.10) | 2.91(0.11) | 3.03(0.08) |
| 2048 | 2.94(0.06) | 3.05(0.04) | 2.90(0.06) | 2.97(0.04) |
| 8192 | 2.92(0.03) | 2.96(0.02) | 2.86(0.03) | 2.95(0.02) |
| 32768 | 2.90(0.02) | 2.94(0.01) | 2.86(0.02) | 2.95(0.01) |

6.3.3. *Asian option*. The dynamics is unchanged (with $d = q = 1$) but now the interest rates are equal $(r = R)$. The terminal condition equals $\Phi(\mathbf{S}) = (\frac{1}{T} \int_0^T S_t \, dt - K)_+$ and we take the following parameters:

| $\mu$ | $\sigma$ | $r$ | $T$ | $S_0$ | $K$ |
|---|---|---|---|---|---|
| 0.06 | 0.2 | 0.1 | 1 | 100 | 100 |

To approximate this path-dependent terminal condition, we take $d' = 2$ and simulate $P_{t_k}^N = (S_{t_k}^N, \frac{1}{k+1} \sum_{i=0}^k S_{t_i}^N)^*$ (see [18]). The results presented in Table 5 are coherent because the price given by the algorithm is not far from the reference price 7.04 given in [18].

As mentionned in [18], the use of $\frac{1}{N+1} \sum_{i=0}^N S_{t_i}^N$ to approximate $\frac{1}{T} \int_0^T S_t \, dt$ is far from being optimal. We can check what happens if we change $P^N$ to better approximate $\frac{1}{T} \int_0^T S_t \, dt$. As proposed in [18], we approximate $\frac{1}{T} \int_0^T S_t \, dt$ by $\frac{1}{N} \sum_{i=0}^{N-1} S_{t_i}^N (1 + \mu \frac{h}{2} + \frac{\sigma}{2} \Delta W_{t_i})$, which leads to $P_{t_k}^N = (S_{t_k}^N, \frac{1}{k} \sum_{i=0}^{k-1} S_{t_i}^N (1 + \mu \frac{h}{2} + \frac{\sigma}{2} \Delta W_{t_i}))^*$ for $k \geq 1$. The results (see Table 6) are much better with this choice of $P^N$. Once more, we observe the coherence of the algorithm which takes in input simulations of $S^N$ under the historical probability ($\mu \neq r$) and corrects the drift to give the risk-neutral price.

TABLE 4
*Results for the calls combination using the basis GP*

| $M$ | $N = 5$ $d_y = 1, d_z = 0$ | $N = 20$ $d_y = 2, d_z = 1$ | $N = 50$ $d_y = 4, d_z = 2$ | $N = 50$ $d_y = 9, d_z = 9$ |
|---|---|---|---|---|
| 128 | 2.87(0.39) | 3.01(0.24) | 3.02(0.22) | 3.49(1.57) |
| 512 | 2.82(0.20) | 2.94(0.12) | 2.97(0.09) | 3.02(0.1) |
| 2048 | 2.78(0.07) | 2.92(0.07) | 2.92(0.04) | 2.97(0.03) |
| 8192 | 2.78(0.05) | 2.92(0.04) | 2.92(0.02) | 2.96(0.01) |
| 32768 | 2.79(0.03) | 2.91(0.02) | 2.91(0.01) | 2.95(0.01) |

TABLE 5
*Results for the Asian option using the basis HC*

| $M$ | $N = 5$ $\delta = 5$ $D = [60, 200]^2$ | $N = 20$ $\delta = 1$ $D = [60, 200]^2$ | $N = 50$ $\delta = 0.5$ $D = [60, 200]^2$ |
|---|---|---|---|
| 128 | 6.33(0.41) | 4.47(3.87) | 3.48(13.08) |
| 512 | 6.65(0.21) | 6.28(0.76) | 5.63(2.37) |
| 2048 | 6.80(0.09) | 6.76(0.24) | 6.48(0.49) |
| 8192 | 6.83(0.04) | 6.95(0.06) | 6.86(0.12) |
| 32768 | 6.83(0.02) | 6.98(0.03) | 6.99(0.04) |

**7. Conclusion.** In this paper we design a new algorithm for the numerical resolution of BSDEs. At each discretization time, it combines a finite number of Picard iterations (3 seems to be relevant) and regressions on function bases. These regressions are evaluated at once with one set of simulated paths, unlike [6], where one needs as many sets of paths as discretization times. We mainly focus on the theoretical justification of this scheme. We prove $\mathbf{L}_2$ estimates and a central limit theorem as the number of simulations goes to infinity. To confirm the accuracy of the method, we only present few convincing tests and we refer to [19] for a more detailed numerical analysis. Even if no related results have been presented here, an extension to reflected BSDEs is straightforward (as in [6]) and allows to deal with American options. At last, we mention that our results prove the convergence of the *Hedged Monte Carlo* method of Bouchaud, Potters and Sestovic [5], which can be expressed in terms of BSDEs with a linear driver.

TABLE 6
*Results for the Asian option, using a better approximation of $\frac{1}{T}\int_0^T S_t\,dt$ and the basis HC ($N = 20$, $\delta = 1$, $D = [60, 200]^2$)*

| $M$ | 2 | 8 | 32 | 128 | 512 | 2048 | 8192 | 32768 |
|---|---|---|---|---|---|---|---|---|
| $\overline{Y}_{t_0}^{N,I,M}$ | 2.26 | 0.90 | 4.49 | 6.68 | 6.15 | 6.88 | 6.99 | 7.02 |
| $\sigma_{t_0}^{N,I,M}$ | 4.08 | 7.80 | 11.27 | 4.64 | 1.11 | 0.21 | 0.07 | 0.02 |

E. GOBET  
ÉCOLE POLYTECHNIQUE  
CENTRE DE MATHÉMATIQUES APPLIQUÉES  
91128 PALAISEAU CEDEX  
FRANCE  
E-MAIL: emmanuel.gobet@polytechnique.fr

J. P. LEMOR  
X. WARIN  
ÉLECTRICITÉ DE FRANCE  
EDF R&D, SITE DE CLAMART  
1 AVENUE DU GÉNÉRAL DE GAULLE  
92141 CLAMART  
FRANCE  
E-MAIL: lemor@cmapx.polytechnique.fr  
xavier.warin@edf.fr